\catcode`\@=11

\magnification=1200
\hsize 16.2truecm      \vsize 22truecm
\hoffset -0,15truecm      \voffset=10mm
\pretolerance=500 \tolerance=1000 \brokenpenalty=5000





\font\tenbb=msbm10
\font\sevenbb=msbm7
\font\fivebb=msbm5
\newfam\bbfam 
\textfont\bbfam=\tenbb  
\scriptfont\bbfam=\sevenbb
\scriptscriptfont\bbfam=\fivebb
\def\bb{\fam\bbfam\tenbb}

\font\tengoti=eufm10 
\font\sevengoti=eufm7 
\font\fivegoti=eufm5 
\newfam\gotifam \textfont\gotifam=\tengoti
\scriptfont\gotifam=\sevengoti
\scriptscriptfont\gotifam=\fivegoti
\def\goth{\fam\gotifam\tengoti}

\font\tensmc=cmcsc10
\def\smc{\tensmc}

\frenchspacing

\catcode`\@=12

\def\hfl#1#2{\smash{\mathop{\hbox to 9mm{\rightarrowfill}}
\limits^{\scriptstyle#1}_{\scriptstyle#2}}}
\def\lfl#1#2{\smash{\mathop{\hbox to 9mm{\leftarrowfill}}
\limits^{\scriptstyle#1}_{\scriptstyle#2}}}

\def\cqfd{\vbox{\hrule 
\hbox to 0.7em{%
\vrule height 1.4ex\hfil
\vrule height 1.4ex 
}\hrule}}


{\hfill September 2017}
\bigskip\bigskip
\centerline {\bf INFINITESIMAL SYMMETRIES FOR IMPLICIT PLANAR WEBS}
\bigskip
\bigskip
\centerline{{\smc Alain H\'enaut}}
\smallskip
\centerline{\it Institut de Math\'ematiques de Bordeaux}

\centerline{Universit\'e de Bordeaux et CNRS, France}

\centerline{\tt Alain.Henaut@math.u-bordeaux.fr}
\bigskip\bigskip
\bigskip
{\bf 1. Introduction}
\smallskip
In the complex setting, we study the geometry of integral curves of an analytic or algebraic differential equation  
$$F(x,y,y')=a_0(x,y)\,.\,(y')^d+\cdots+a_{d-1}(x,y)\,.\,y'+a_d(x,y)=0$$ of the first order with degree $d\geq 1$. In $({\bb C}^2,0)$ with the germ language, $F\in {\cal O}[p]$ where $\displaystyle p={dy\over dx}$ and
${\cal O}:={\bb C}\{x,y\}$ is the ring of convergent power series in two
variables or globally in ${\bb P}^2:={\bb P}^2({\bb C})$ by using the changes of standard affine charts, $F\in {\bb C}[x,y,p]$. To be more precise, we suppose that the $p$-resultant of $F$ satisfies $R_F:=\hbox{Result}(F,\partial_p(F))=(-1)^{d(d-1)\over 2}a_0\,.\,\Delta\not=0$, where $\Delta\in
{\cal O}$ is the $p$-discriminant of $F$ and, to avoid parasit
solutions, that all the $a_i$ are
relatively prime.
 
If $R_F(0)\not=0$, we have $\displaystyle F=a_0\prod_{i=1}^d(p-p_i)$ and we get a (germ of a) {\it nonsingular $d$-web}
${\cal W}(d)={\cal W}(F_1,\ldots,F_d)$ in $({\bb
C}^2,0)$, that is $d$ foliations by curves in general position. Its leaves are level sets of elements $F_i\in {\cal O}$ where
$F_i(0)=0$ and $X_i(F_i)=0$ with analytic vector fields
$X_i=\partial_x+p_i\,\partial_y\in \Theta$ or analytic 1-forms $dy-p_idx\in \Omega^1$ corresponding  to the different
slopes $p_i=p_i(x,y)$ of $F(x,y,y')=0$. Conversely and up to a linear transformation to avoid the ``vertical'' slopes, every ${\cal W}(F_1,\ldots,F_d)$ gives rise to a differential equation $F(x,y,y')=0$ as above by using the product of the $\partial_y(F_i)\,.\,p+\partial_x(F_i)$.

Planar web geometry studies the previous configurations ${\cal W}(d)$, up to local isomorphisms ({\it cf.} for example the classic reference [BB-1938] or [PP-2015] which deals with nonsingular webs of codimension one in general). Here $F$ and $g\,.\,F$ represent the same web ${\cal W}(d)$ for $g\in {\cal O}^*$ and the main goal is to study the geometry of the class of the differential equations defined by $F(x,y,y')=0$ and their singularities by using the generic ${\cal W}(d)$ associated. In the implicit setting the previous integral curves, as the roots of a polynomial equation, will be investigated in the {\it mutual relations} of all. This subject can also be viewed as a part of the {\it qualitative study} of differential equations, especially the binary differential equations of degree $d$.
\smallskip
In the following we deal with infinitesimal symmetries or to be short symmetries for such a ${\cal W}(d)$. With the previous notations, a {\it symmetry} is a vector field $X=\alpha\partial_x+\beta\partial_y$ such the local flow $(x,y\,;t)\longmapsto \exp(tX)(x,y)$ generated by $X$ preserves all the leaves of ${\cal W}(d)$. This means generically that $X_i\big(X(F_i)\big)=0$ or equivalently ${\cal
L}_X(\omega_i)\wedge\omega_i=0$ for $1\leq i\leq d$ where the 1-forms $\omega_i$ define ${\cal
W}(d)$. Here ${\cal
L}_X=i_X\circ d+d\circ i_X:\Omega^\bullet\longrightarrow \Omega^\bullet$ is the Lie derivative associated with $X$ and $i_X:\Omega^\bullet\longrightarrow \Omega^{\bullet -1}$ is its interior product. The symmetries $\cal S$ of a ${\cal W}(d)$ form a sheaf such that each fiber, equipped with the Lie bracket, is a Lie algebra. Elements in $\cal S$ are also called Lie symmetries or infinitesimal automorphisms of the ${\cal W}(d)$ at stake.  
\smallskip
For special linear webs in ${\bb P}^2$ associated by duality with plane reduced algebraic curves, the study of the associated symmetries is directly related to W-curves first investigated by Felix Klein and Sophus Lie in [KL-1871] as an important step in the emergence of the theory of Lie groups. As an illustration of the equivalence problem, \'Elie Cartan in [C-1908] studies the generic case of a planar $3$-web ${\cal W}(x,y,F_3)$ and its associated $\cal S$. He pointed out that $\dim {\cal S}$, as $\bb C$-vector space, is $0$, $1$ or $3$. In particular the basic hexagonal $3$-web ${\cal W}(x,y,x+y)$ corresponds to $\dim {\cal S}=3$. Infinitesimal symmetries do not appear in the pioneer work on web geometry by Wilhelm Blaschke and his coworkers in the 1930's ({\it cf.} [BB-1938]). It is also the case in the basic revival of interest in web geometry due to Shiing-Shen Chern and Phillip A. Griffiths which essentially focus on abelian relations for webs and algebraization problems ({\it cf.} [Ch-1982] and references therein). Some results on these symmetries for a $3$-web of codimension $r$ in a $2r$-dimensional manifold can be found in the contributions of Maks A. Akivis with Alexander M. Shelekhov [AS-1992] or Vladislav V. Goldberg [AG-2000] (this survey contains a detailed report of the publications on this subject in russian by N. V. Gvozdovich in the 1980's). More recently David Mar\'\i n, Jorge Vit\'orio Pereira and Luc Pirio use these symmetries coupled with plane algebraic W-curves in [MPP-2006] to provide a large class of remarkable planar webs ${\cal E}(d)$ for $d\geq 5$. These webs are not polylogarithmic as Bol's ${\cal B}(5)$ which for example have no symmetry in the previous sense. Sergey I. Agafonov by precising results of Eugene V. Ferapontov also studies $\cal S$ for particular planar $3$-webs which are related to the geometry of Frobenius manifolds introduced by Boris Dubrovin through WDVV-equations ({\it cf.} [Ag-2012] for example). Some part of these previous contributions will be specified or even revisited in the following.
\smallskip
It is proved below that the symmetries $\cal S$ of a ${\cal W}(d)$ presented by $F$ with $d\geq 3$ form a {\it local system} with rank $0$, $1$ or $3$, outside its discriminant locus $|\Delta|$ defined by the reduced divisor associated with $\Delta$. For $d=3$, the local system  $\cal S$ is incarnated as horizontal sections of a {\it connection} $(E,\nabla)$ with rank $3$ which is meromorphic on $|\Delta|$ and integrable (or flat depending on the terminology) if, and only if, ${\cal W}(3)$ is hexagonal. A complete discussion on $\dim {\cal S}$ is given with {\it effective methods} depending only on the coefficients of $F$ in case $d=3$. Generalization in case $d\geq 3$ are also presented. The implicit viewpoint where no leaf of the web at stake is prefered plays its role, in the nonsingular case and through the singularities as well. New invariants of the class of the differential equations defined by $F(x,y,y')=0$ which extend the classic curvature of a planar $3$-web are introduced by connection methods already used in [H\'e-2004]. Links between $\cal S$ and rank problems related to abelian relations of such a ${\cal W}(d)$ are also studied, especially with regards to remarkable planar webs. A generalization for a planar web of Lie's equivalence for a 1-form between existence of an integrating factor and a transverse symmetry is presented. Detailed examples, discussions on singularities of $\cal S$ including regularity aspects and the presence of weighted Euler symmetry and some perspectives complete the present article. 
\smallskip
The author would like to thank Daniel Lehmann for general comments on the first draft of the present text and also David Mar\'\i n who moreover make available an illuminating family of hexagonal webs in ${\bb P}^2$ which in particular is not always regular singular ({\it cf.} Example 6, below).   
\bigskip
{\bf 2. Description of $\cal S$ as a local system and Lie's integrating factor for an implicit planar web}
\smallskip 
In this section, we suppose that $R_F(0)\not=0$ where $R_F=\hbox{Result}(F,\partial_p(F))$.
 
By definition and from a calculus using the forms $\omega_i=dy-p_idx$ or the vector fields $X_i=\partial_x+p_i\partial_y$ which define the nonsingular
${\cal W}(d)$, we get that $X=\alpha\partial_x+\beta\partial_y\in \Theta$ is a symmetry for ${\cal W}(d)$ that is $X\in {\cal S}$ if, and only if, $(\alpha,\beta)$ is an analytic solution of the following (homogeneous) linear differential system with $d$ equations:
$$({\cal S})\quad-\partial_x(\beta)+\big(\partial_x(\alpha)-\partial_y(\beta)\big)\,.\,p_i+\partial_y(\alpha)\,.\,p_i^2+\partial_x(p_i)\,\alpha+\partial_y(p_i)\,\beta=0\ \ \,\hbox{for}\  1\leq i\leq d.$$
\smallskip
For $d=1$ or $d=2$, the $\bb C$-vector space $\cal S$ is not finite dimensional. Indeed, in these cases we may suppose $p_1=0$, or $p_1=0$ with $p_2=1$. Hence, every $X=\beta(y)\partial_y$ with $\beta\in {\bb C}\{y\}$ is a solution in the first case while every $X=\beta(y)(\partial_x+\partial_y)$ is a solution in the second case. For $d\geq 3$, the system $({\cal S})$ is equivalent to a system which looks like to the linear differential system ${\cal M}(4)$ with 3 equations associated with the presentation of a $4$-planar web in [Hé-2004] with an additional part formed by $d-3$ linear equations. This linear part can be viewed as compatibility conditions associated with the differential system at stake. 

More precisely for $d\geq 3$, that we suppose from now on, we have the following equivalence: 
\smallskip
$X=\alpha\partial_x+\beta\partial_y\in {\cal S}$ {\it if, and only if, $(\alpha,\beta)$ is an analytic solution of the following linear differential system}:
$$({\cal S})\quad\left\{\matrix{
-\,\partial_x(\beta)&+&g_d\alpha+h_d\beta&=&0\hfill\cr
\partial_x(\alpha)-\partial_y(\beta)&+&g_{d-1}\alpha+h_{d-1}\beta&=&0\hfill\cr
\partial_y(\alpha)&+&g_{d-2}\alpha+h_{d-2}\beta&=&0\hfill\cr
&&g_{d-3}\alpha+h_{d-3}\beta&=&0\hfill\cr
&&&\ \vdots\hfill\cr
&&g_{1}\alpha+h_{1}\beta&=&0.\hfill\cr
}\right.$$

The previous $(d\times 2)$-matrix $(S_{ij}):=\pmatrix{
g_d&h_d\cr
\vdots&\vdots\cr
g_1&h_1\cr
 }$ is called the {\it matrix of symmetries} of the presentation $F$ of ${\cal W}(d)$. It comes from the two polynomials 
$$G:=g_1\,.\,p^{d-1}+\cdots +g_{d-1}\,.\,p+g_d\quad \hbox{and}\ \ H:=h_1\,.\,p^{d-1}+\cdots +h_{d-1}\,.\,p+h_d$$ which are the unique elements in ${\cal O}[p]$ with $\deg G\leq d-1$ and $\deg H\leq d-1$ such that $G(x,y,p_i)=\partial_x(p_i)$ and $H(x,y,p_i)=\partial_y(p_i)$ for $1\leq i\leq d$. Indeed, the initial description of $({\cal S})$ can be also written, as the following system: 
$$\eqalign{
-\partial_x(\beta)
+g_d\alpha+h_d\beta\,\ \ \ &\cr 
+\big(\partial_x(\alpha)-\partial_y(\beta)+g_{d-1}\alpha+h_{d-1}\beta\big)\,.\,&p_i\cr
+\big(\partial_y(\alpha)+g_{d-2}\alpha+h_{d-2}\beta\big)\,.\,&p_i^2\cr
+\big(g_{d-3}\alpha+h_{d-3}\beta\big)\,.\,&p_i^3
\cr
+\cdots
+\big(g_{1}\alpha+h_{1}\beta\big)\,.\,&p_i^{d-1}=0\cr
}$$
for $1\leq i\leq d$. Then, using the Vandermonde $d$-determinant obtained from the distinct $p_i$, we get the equivalence below.
\smallskip
The $p$-resultant $R_F=\hbox{Result}(F,\partial_p(F))\in {\cal O}$ is classically interpreted as a $(2d-1)$-determinant. This proves that the previous $G$ and $H$ can be also viewed as part of the unique ordered pair of polynomials $(V,G)$ and $(W, H)$ in ${\cal O}[p]_{<d-1}\times {\cal O}[p]_{<d}$ such that
$$V\,.\,F-G\,.\,\partial_p(F)=\partial_x(F)\quad\hbox{and}\quad W\,.\,F-H\,.\,\partial_p(F)=\partial_y(F)$$since $F(x,y,p_i)=0$ for $1\leq i\leq d$ by definition. Here ${\cal O}[p]_{<\ell}$ denotes the elements in  ${\cal O}[p]$ with a $p$-degree strictly less than $\ell$. 

Therefore to be concrete and for $d=3$, we have the following Cramer-Sylvester systems: 
$$\pmatrix{
a_0&0&3a_0&0&0\cr
a_1&a_0&2a_1&3a_0&0\cr
a_2&a_1&a_2&2a_1&3a_0\cr
a_3&a_2&0&a_2&2a_1\cr
0&a_3&0&0&a_2\cr
}\pmatrix{
*\cr
*\cr
-g_1 \cr
-g_2 \cr
-g_3 \cr
}=\pmatrix{
0\cr
\partial_x(a_0)\cr
\partial_x(a_1)\cr
\partial_x(a_2)\cr
\partial_x(a_3)\cr
 }$$ and $$\pmatrix{
a_0&0&3a_0&0&0\cr
a_1&a_0&2a_1&3a_0&0\cr
a_2&a_1&a_2&2a_1&3a_0\cr
a_3&a_2&0&a_2&2a_1\cr
0&a_3&0&0&a_2\cr
}\pmatrix{
*\cr
*\cr
-h_1 \cr
-h_2 \cr
-h_3 \cr
}=\pmatrix{
0\cr
\partial_y(a_0)\cr
\partial_y(a_1)\cr
\partial_y(a_2)\cr
\partial_y(a_3)\cr
 }$$
where $G=g_1\,.\,p^2+g_2\,.\,p+g_3$ and $H=h_1\,.\,p^2+h_2\,.\,p+h_3$. For $d=4$, we get  
$$\pmatrix{
a_0&0&0&4a_0&0&0&0\cr
a_1&a_0&0&3a_1&4a_0&0&0\cr
a_2&a_1&a_0&2a_2&3a_1&4a_0&0\cr
a_3&a_2&a_1&a_3&2a_2&3a_1&4a_0\cr
a_4&a_3&a_2&0&a_3&2a_2&3a_1\cr
0&a_4&a_3&0&0&a_3&2a_2\cr
0&0&a_4&0&0&0&a_3\cr
}\pmatrix{
*\cr
*\cr
*\cr
-g_1 \cr
-g_2 \cr
-g_3 \cr
-g_4 \cr
}=\pmatrix{
0\cr
0\cr
\partial_x(a_0)\cr
\partial_x(a_1)\cr
\partial_x(a_2)\cr
\partial_x(a_3)\cr
\partial_x(a_4)\cr
 }$$
and the analogue for $H$, and so on for $d\geq 5$. 

It can be noted that the coefficients $S_{ij}$ of the matrix of symmetries depend only on $a_i$, $\partial_x(a_i)$ and $\partial_y(a_i)$. From  the singularity viewpoint these are in fact in ${\cal O}[1/\delta]$, that is with poles on the {\it discriminant locus} $|\Delta|$ of ${\cal W}(d)$ defined by the reduced divisor $\displaystyle\delta:=\prod_q \Delta_q$ associated with the irreducible decomposition $\displaystyle \Delta=u\,.\,\prod_q \Delta_q^{m_q}$ of the 
$p$-discriminant of $F$ where $u\in {\cal O}^*$ and every 
$\Delta_q$ defines an irreducible analytic germ in $({\bb C}^2,0)$ with
$m_q\geq 1$.
\medskip
\noindent
{\bf Remark 1.} \  The geometric definition of a symmetry $X\in \cal S$ through its generated flow does not depend on the local coordinates where ${\cal W}(d)$ is presented. The polynomials $G$ and $H$ depend on the coordinates, but the presentations $F$ and $g\,.\,F$ of the same ${\cal W}(d)$ for $g\in {\cal O}^*$ give the same $G$ and $H$. These properties can be viewed also as a direct consequence of the Lagrange interpolation formula. In other words, the polynomials $G$ and $H$ are {\it invariants} of the class of the differential equations defined by $F(x,y,y')=0$.   
\medskip
\noindent
{\bf Remark 2} ({\it Linear webs and web in ${\bb P}^2$ associated by duality with a reduced algebraic curve}){\bf.}\ If  ${\cal W}(d)={\cal L}(d)$ is a {\it linear $d$-web}, that is  all its leaves are straight lines, we have equivalently
$$X_i(p_i)=\partial_x(p_i)+p_i\partial_y(p_i)=0\ \ \hbox{for} \ \ 1\leq i\leq d$$
by using the local flow generated by the vector field $X_i$. Therefore in this case we obtain
$$a_0 g_d=a_{d}h_1\ ,\ a_0(g_{d-1}+h_d)=a_{d-1}h_1\ ,\ \ldots\ , \ a_0(g_{1}+h_2)=a_1h_1.$$
Indeed, if $h_1=0$, these equalities are coming from the Vandermonde $d$-determinant already used by definition of $G$ and $H$. If $h_1\not=0$, that is $\deg H=d-1$, then the polynomial $G+p\,.\,H\in {\cal O}[p]$ has the same roots $p_i$ that the presentation $F$ of ${\cal L}(d)$ and the previous equalities follows. Let $C\subset \check{\bb P}^2:=G(1,{\bb P}^2)$ be a reduced algebraic curve with degree $d$ in the space of lines of ${\bb P}^2$. By duality, we get a special linear $d$-web ${\cal L}_C(d)\subset {\bb P}^2$ generically nonsingular, called {\it the $d$-web associated with} $C$. It is locally presented by 
$$F(x,y,p)=P(y-px,p)$$ 
if $P(q,p)=0$ is an affine equation of $C$. Here $F(x,y,y')=0$ corresponds essentially to $d$ classic Clairaut's equations $y=xy'+f_i(y')$. If $C$ contains no lines, the leaves of ${\cal L}_C(d)$ are generically the tangents of the {\it dual curve} $\check C\subset {\bb P}^2$ of $C$. This one is locally defined by $\check P(x,y)=0$ where $\check P$ is a factor of the $p$-resultant $R_F$, and the others factors are products of linear forms related to the singular points of $C$. If otherwise $C$ contains lines, then the corresponding points in ${\bb P}^2$ give rise to pencils of lines for ${\cal L}_C(d)$ through these points.     
\medskip
\noindent
{\bf Example 1.} {\it Parallel planar $d$-web}. 
\smallskip
We consider a particular $d$-web ${\cal SP}(d)\subset {\bb P}^2$ generated by the {\it special pencils} of lines passing through $d$ distinct points which belong to a line $L\subset {\bb P}^2$. By duality ${\cal SP}(d)={\cal L}_C(d)$ where $C\subset \check{\bb P}^2$ is the algebraic curve given by the union of $d$ distincts lines passing through the point $\check L\in\check{\bb P}^2$, that is a central arrangement of $d$ lines. With notations of Remark 2, we may suppose that $C$ is given by $\displaystyle P(q,p)=\prod_{i=1}^d(p-p_i)=0$ where the $p_i$ are $d$ distinct complex numbers and $\check L=[1,0,0]$ or $L=\{[X_0,X_1,X_2]\,;\,X_0=0\}$. From the previous definitions, the following properties are equivalent:
\smallskip
$i)$ {\it ${\cal W}(d)={\cal SP}(d)$ is a parallel planar $d$-web, that is presented by $F\in {\bb C}[p]$;} 
\smallskip
$ii)$ {\it The slopes $p_i$ of ${\cal W}(d)$ are $d$ distinct constants;}
\smallskip
$iii)$ {\it The $(d\times 2)$-matrix of symmetries of $F$ is $(S_{ij})=0$, that is $G=H=0$.}
\smallskip
\noindent  
In this case the system $({\cal S})$ is reduced to the following:
$$({\cal S})\quad\left\{\matrix{
-\,\partial_x(\beta)&=&0\hfill\cr
\partial_x(\alpha)-\partial_y(\beta)&=&0\hfill\cr
\partial_y(\alpha)&=&0.\hfill\cr
}\right.$$
Therefore we obtain ${\cal S}=\{\partial_x,\,\partial_y,\,x\partial_x+y\partial_y\}$ for such a ${\cal SP}(d)$. 
\medskip
Since in this part $R_F(0)\not=0$, we have a nonsingular surface $S$ defined by $F$ with its usual de Rham complex $(\Omega^{\bullet}_S,d)$ where 
$\Omega^{\bullet}_S=\Omega^{\bullet}_{{\bb C}^3}/(dF\wedge
\Omega^{\bullet-1}_{{\bb C}^3},F\,\Omega^{\bullet}_{{\bb C}^3})$. Moreover the first projection induces a $d$-covering map $\pi:S\longrightarrow ({\bb C}^2,0)$ with local branches $\pi_i(x,y)=(x,y,p_i(x,y))$. 
\smallskip
The following result summarizes the previous observations. It also generalizes classical Lie's link for a 1-form between existence of an integrating factor and a transverse symmetry, to an implicit planar $d$-web. The classic link is precisely the equivalence $\dot v )\Longleftrightarrow ii)$ proved below for a given $\omega_i$.
\medskip 
\noindent
{\bf Theorem 1.}\ \ {\it With the previous notations, in particular $\omega_i=dy-p_idx$ and $X=\alpha\partial_x+\beta\partial_y\in \Theta$, the following conditions are equivalent for a  ${\cal W}(d)$ with $d\geq 3$:
\smallskip
$i)$ $X\in {\cal S}$, that is $X$ is a symetry for ${\cal W}(d)$;
\smallskip
$ii)$ ${\cal
L}_X(\omega_i)\wedge\omega_i=0$ for $1\leq i\leq d$;
\smallskip
$iii)$ $(\alpha,\beta)$ is an analytic solution of the linear differential system $({\cal S})$.
\smallskip\noindent
Moreover if $i_X(\omega_i)\not= 0$ for $1\leq i\leq d$, that is $X$ is transverse to ${\cal W}(d)$, then the previous conditions are equivalent to the following:
\smallskip 
$iv)$  The meromorphic $1$-form $\omega:=\displaystyle {dy-pdx\over {\beta-\alpha p}}$ on $S$ is closed, that is
$\displaystyle {1\over  {\beta-\alpha p}}$ is an integrating factor for $dy-pdx$ on the surface $S$;
\smallskip
$\dot v)$ $\displaystyle {1\over  i_X(\omega_i)}$ is an integrating factor for $\omega_i$ for $1\leq i\leq d$}.
\medskip
\noindent
{\it Proof}. The previous observations prove that $i)\Longleftrightarrow ii)\Longleftrightarrow iii)$.

\noindent
$iii)\Longrightarrow iv)$. By hypothesis, $i_X(\omega_i)=\omega_i(X)=\beta-\alpha p_i\not= 0$ on $S$. In $\Omega^2_S$, we have $\displaystyle {dx\wedge dy\over \partial_p(F)}= {dy\wedge dp\over \partial_x(F)}= {dp\wedge dx\over \partial_y(F)}$ since $\partial_x(F)dx+\partial_y(F)dy+\partial_p(F)dp=0$ in $\Omega^1_S$. Hence a calculus gives
$$\eqalign{
d\omega&=\Big(p^2\partial_y(\alpha)+p\big(\partial_x(\alpha)-\partial_y(\beta)\big) -\partial_x(\beta)-\alpha\,{\partial_x(F)\over \partial_p(F)}-\beta\,{\partial_y(F)\over \partial_p(F)}\Big){dx\wedge dy\over (\beta-\alpha p)^2} \cr
&=\Big(p^2\partial_y(\alpha)+p\big(\partial_x(\alpha)-\partial_y(\beta)\big) -\partial_x(\beta)+\alpha\cdot G +\beta\cdot H\Big){dx\wedge dy\over (\beta-\alpha p)^2}=0 \cr
}$$ on the surface $S$ by definition of $G$ and $H$ through the expressions above.

\noindent
$iv)\Longrightarrow \dot v)$. For $1\leq i\leq d$, we have $\displaystyle \pi_i^*(\omega)={\omega_i\over i_X(\omega_i)}$ by definition and this form is closed by hypothesis since $d$ and $\pi_i^*$ commute. 

\noindent
$\dot v)\Longrightarrow ii)$. For $1\leq i\leq d$, we have $0=i_X(\omega_i\wedge d\omega_i)=i_X(\omega_i)d\omega_i-\omega_i\wedge i_X(d\omega_i)$ and by hypothesis
$\,\displaystyle d({\omega_i\over i_X(\omega_i)})={d\omega_i\over i_X(\omega_i)}-{d(i_X(\omega_i))\over i_X(\omega_i)^2}\wedge \omega_i=0$. Therefore, we get  
$${\cal L}_X(\omega_i)\wedge \omega_i=\big(i_X(d\omega_i)+d(i_X(\omega_i)\big)\wedge \omega_i=-i_X(\omega_i)d\omega_i+d(i_X(\omega_i))\wedge \omega_i=0.\ \cqfd$$ 
\smallskip
Discussions on $\dim {\cal S}$ will be given in the next paragraphs. However, we begin with the following result:
\medskip 
\noindent
{\bf Proposition 1.}\ \ {\it Let ${\cal W}(d)$ be an implicit planar $d$-web presented by $F$ with a discriminant locus $|\Delta|$. Then the sheaf $\cal S$ of symmetries of ${\cal W}(d)$ is a local system outside  $|\Delta|$, that is a locally constant sheaf of ${\bb C}$-vector spaces with finite dimensional fibers $\dim {\cal S}$. Moreover $\dim {\cal S}$ takes values in $\{0,1,3\}$}.
\medskip
\noindent
{\it Proof.} With usual notations and basic results in $\cal D$-modules approach where $\cal D$ is the ring of linear differential operators with coefficients in $\cal O$ ({\it cf.} for example [GM-1993]) and after transposition, the symbol matrix of the differential system in the $3$-upper part of $({\cal S})$ is
$$\pmatrix{
0&\xi&\eta\cr
-\xi&-\eta&0\cr
}.$$ Here the left $\cal D$-module at stake is the cokernel of $^t\rho :{\cal D}^3\longrightarrow {\cal D}^2$ where  $^t\rho(Q_1,Q_2,Q_3)=(Q_1g_d+Q_2(\partial_x+g_{d-1})+Q_3(\partial_y+g_{d-2}),Q_1(-\partial_x+h_d)+Q_2(-\partial_y+h_{d-1})+Q_3h_{d-2})$. Therefore the 0-th Fitting ideal associated, that is the ideal of $2\times2$-minors generated by $^t\rho$ in $\hbox{gr}\,{\cal D}={\cal O}[\xi,\eta]$ is $(\xi^2,\xi\eta,\eta^2)$. Hence as left $\cal D$-modules, we have $\hbox{Coker}\,^t\rho={\cal O}^m$ with $m\leq 3=\hbox{mult}\,{\cal O}[\xi,\eta]/(\xi^2,\xi\eta,\eta^2)$. Therefore ${\cal S}$ or the analytic solutions of the system $({\cal S})$ is a local system outside the zero locus of $\Delta$ and $\dim {\cal S}$ is bounded by 3. If ${\cal S}$ has dimension 2, we may suppose ${\cal S}=\{\partial_x, \partial_y\}$ with ${\cal L}_{\partial_x}(\omega_i)\wedge \omega_i=0= {\cal L}_{\partial_y}(\omega_i)\wedge \omega_i$ for $1\leq i \leq d$, hence we obtain that $G=H=0$ by definition since the previous conditions are equivalent with $\partial_x(p_i)=0=\partial_y(p_i)$ for $1\leq i \leq d$. In other words, ${\cal W}(d)$ is parallelizable ({\it cf.} Example 1). Therefore, we get also $x\partial_x+y\partial_y\in {\cal S}$ by using the system $({\cal S})$, which is a contradiction. \cqfd
\medskip 
Symmetries $X$ for a singular $d$-web ${\cal W}(d)$ in $({\bb C}^2,0)$ are multivalued in general if there exist and the algebraic analysis of these merits to be undertaken, especially regularity and monodromy aspects, as some examples below show. Those complex analytic or holomorphic, that is with $X\in \Theta$ are interesting in regard to the singularities of ${\cal W}(d)$ ({\it cf.} [Ag-2015] for a particular case) and their different types of indices attached. This vein will be subsequently explored and related to Kyoji Saito's logarithmic vector fields, that is the free $\cal O$-module with rank $2$ defined by $\hbox{Der}(\log |\Delta|)=\{X\in \Theta\,;\,X(\delta)\in (\delta)\}$ in [S-1980]. In particular we mention the general invariance property which follows.
\medskip 
\noindent
{\bf Proposition 2.}\ \ {\it Let ${\cal W}(d)$ be an implicit planar $d$-web presented by $F=a_0\,.\,p^d+a_1\,.\,p^{d-1}+\cdots+a_d$ which admits a symmetry $X=\alpha\partial_x+\beta\partial_y$. Then the $p$-discriminant $\displaystyle \Delta=u\,.\,\prod_q \Delta_q^{m_q}$ of $F$ is invariant by $X$. More precisely, we have
$$X(\Delta)=\lambda_X\,.\,\Delta$$
with $\displaystyle \lambda_X=(d-1)\,.\,\Big({2X(a_0)\over a_0}-d\,.\,\big(\partial_x(\alpha)-\partial_y(\beta)\big)+2\partial_y(\alpha)\,{a_1\over a_0}\Big)$. In particular outside $|\Delta|$ with $u\equiv 1$, there exist $\lambda_{X,q}\in {\cal O}$ such that $X(\Delta_q)=\lambda_{X,q}\,.\,\Delta_q$, that is every irreducible component $\Delta_q$ of $\Delta$ is invariant by $X$ with $\displaystyle \lambda_X=\sum _qm_q\,.\,\lambda_{X,q}$.}
\medskip
\noindent
{\it Proof.} With the previous notations, we have classically $\displaystyle \Delta=a_0^{2d-2}\,.\prod_{1\leq i<j\leq d}(p_i-p_j)^2$. Since $X\in {\cal S}$, we have $X(p_i-p_j)=-\big(\partial_x(\alpha)-\partial_y(\beta)\big)(p_i-p_j)-\partial_y(\alpha)(p_i^2-p_j^2)$ from the initial presentation of the linear differential system $({\cal S})$. Then the invariance formula follows by applying $X$ on the product above and using the classical relation on the sum of the roots $p_i$ of $F$ with its coefficients. In particular for $u\equiv 1$, we obtain the equality
$$\sum_qm_q\,.\,X(\Delta_q)\Delta_1\ldots\Delta_{q-1}\,\widehat{\Delta_q}\,\Delta_{q+1}\ldots\Delta_r=\lambda_X\,.\,\prod_q \Delta_q$$
 outside $|\Delta|$. Hence for example $X(\Delta_1)\Delta_2\ldots \Delta_r$ divides $\Delta_1$ and, by successively using the $\Delta_q$ for $q\not=1$, $X(\Delta_1)$ divides $\Delta_1$ since the $\Delta_q$ are pairwise coprime. So we get the $\lambda_{X,q}$ as claimed. \cqfd
\bigskip
{\bf 3. Basic results on $\cal S$ in case $d=3$}
\smallskip
We first recall some results for an implicit planar 3-web ${\cal W}(3)$. A relation $\xi_1(F_1)dF_1+\xi_2(F_2)dF_2+\xi_3(F_3)dF_3=0$ with $\xi_i\in {\bb C}\{t\}$ is called an {\it abelian relation} of a nonsingular $3$-web
${\cal W}(3)={\cal W}(F_1,F_2,F_3)$. These relations, between the normals of ${\cal W}(3)$, viewed as special 3-uple $\big(\xi_1(F_1),\xi_2(F_2),\xi_3(F_3)\big)\in {\cal O}^3$ form a local system ${\cal A}(3)$ such its rank, noted $\hbox{rank}\,{\cal W}(3)$, is equal to 0 or 1 and does not depend on the choice of the $F_i$. The main invariant of such a nonsingular $3$-web
${\cal W}(3)$, up to the equivalence $\sim$ induces as pullback by an analytic isomorphism of $({\bb
C}^2,0)$, is its {\it curvature} 2-form $k$ also called its Blaschke-Dubourdieu 
curvature. The birth certificate of planar web geometry may be summarized as follows. ({\it cf.} for example [BB-1938] or [PP-2015]).
\medskip
\noindent
{\bf Theorem 0.}\ {\it For a
nonsingular $3$-web ${\cal W}(3)$ in $({\bb C}^2,0)$, the following assertions are
equivalent:} 
\smallskip
$i)$ $\hbox{rank}\,{\cal W}(3)=1$;
\smallskip
$ii)$ ${\cal
W}(3)\sim {\cal
W}(x,y,x+y)$;
\smallskip
$iii)$ {\it ${\cal
W}(3)$ is hexagonal} (or satisfies Thomsen's closure, 1927){\it;} 
\smallskip
$iv) $ $k=0$. 
\medskip
For a ${\cal W}(3)$ presented by $F=a_0\,.\,p^3+a_1\,.\,p^2+a_2\,.\,p+a_3$ there exists an
explicit 
{\it meromorphic $1$-form $\gamma=Adx+Bdy$}, with poles on $|\Delta|$, constructs only from the coefficients
$a_i$, $\partial_x(a_i)$ and  $\partial_y(a_i)$ such that $k=d\gamma$. Description and basic properties of this 1-form $\gamma$ defined in  [Hé-2000] ({\it cf.} also [Hé-2004]) are recalled below.
\smallskip
With the previous notations if $R_F(0)\not=0$, the 1-form $\gamma\in \Omega^1$ is defined by using an explicit description of $d\omega$  on the nonsingular surface $S$ defined by $F$ for $\displaystyle \omega={r_xdy-r_ydy\over \partial_p(F)}\cdot$ Indeed for the normalized contact 1-form $\displaystyle\nu:={{dy-pdx}\over \partial_p(F)}$ on $S$, we get  
$\displaystyle d{\nu}=(A+B\,.\,p)\,{dx\wedge dy\over \partial_p(F)}$ where $A:=-u_3+l_2$ and $B:=-u_2+2l_1$ 
are coming from the
following Cramer-Sylvester system: 
$$\pmatrix{
a_0&0&3a_0&0&0\cr
a_1&a_0&2a_1&3a_0&0\cr
a_2&a_1&a_2&2a_1&3a_0\cr
a_3&a_2&0&a_2&2a_1\cr
0&a_3&0&0&a_2\cr
}\pmatrix{
u_2\cr
u_3\cr
-l_1 \cr
-l_2 \cr
-l_3 \cr
}=\pmatrix{
\partial_y(a_0)\cr
\partial_x(a_0)+\partial_y(a_1)\cr
\partial_x(a_1)+\partial_y(a_2)\cr
\partial_x(a_2)+\partial_y(a_3)\cr
\partial_x(a_3)\cr
 }.$$
 The 1-forms $\displaystyle \nu_i:=\pi_i^*(\nu)={{dy-p_idx}\over \partial_p(F)(x,y,p_i)}$ define the 3-web ${\cal W}(3)$, therefore $k=d\gamma$ is by definition its curvature since $\nu_1+\nu_2+\nu_3=0$ from the Lagrange interpolation formula and $d\nu_i=\gamma\wedge \nu_i$ for $1\leq i\leq 3$. 
Moreover, we have the relation  
$$(*)\quad\quad\quad\partial_x(F)+p\,\partial_y(F)+P_0\,\partial_p(F)=\big(\partial_p(P_0)-A-B\,.\,p\big)\,F$$
where  $P_0:=l_1\,.\,p^2+l_2\,.\,p+l_3$ is the unique element in ${\cal O}[p]$  with $\deg P_0\leq 2$ such that $P_0(x,y,p_i)=X_i(p_i)$ for $1\leq i\leq 3$. The polynomial $P_0$ depends {\it only} on the class of $F(x,y,y')=0$ and the leaves of its associated ${\cal W}(3)$ verify $y''=P_0(x,y,y')$. 
\smallskip
The coefficients $l_j$ and $u_i$ defined from the Cramer-Sylvester
system above belong to ${\cal
O}[1/\hbox{red}\,R_F]$, that is have poles on the $p$-resultant of $F$. However the poles of $A$ and $B$ are only on the discriminant locus $|\Delta|$ with moreover the following detailed expressions:
$$\eqalign{
R_F:&=\hbox{Result}(F,\partial_p(F))=-a_0\,.\,\Delta\cr
\Delta&=-4a_0a_2^3+18a_0a_1a_2a_3-27a_0^2a_3^2+a_1^2a_2^2-4a_1^3a_3\cr
\Delta\,.\,A&=a_3(2a_2^2-6a_1a_3)\,.\,\partial_y(a_0)+a_3(9a_0a_3-a_1a_2)\,.\,(\partial_x(a_0)+\partial_y(a_1))\cr
&\quad\quad+a_3(2a_1^2-6a_0a_2)\,.\,(\partial_x(a_1)+\partial_y(a_2))\cr
&\quad\quad+(4a_0a_2^2-a_1^2a_2-3a_0a_1a_3)\,.\,(\partial_x(a_2)+\partial_y(a_3))\cr
&\quad\quad+(18a_0^2a_3-8a_0a_1a_2+2a_1^3)\,.\,\partial_x(a_3)\cr
\Delta\,.\,B&=(18a_0a_3^2-8a_1a_2a_3+2a_2^3)\,.\,\partial_y(a_0)\cr
&\quad\quad+(4a_1^2a_3-a_1a_2^2-3a_0a_2a_3)\,.\,(\partial_x(a_0)+\partial_y(a_1))\cr
&\quad\quad+a_0(2a_2^2-6a_1a_3)\,.\,(\partial_x(a_1)+\partial_y(a_2))\cr
&\quad\quad+a_0(9a_0a_3-a_1a_2)\,.\,(\partial_x(a_2)+\partial_y(a_3))+a_0(2a_1^2-6a_0a_2)\,.\,\partial_x(a_3).\cr
}$$

It can be noted that for $g\in {\cal O}^*$, the presentation $g\,.\,F$ gives a 1-form $^g\gamma$ where  
$$\displaystyle
^g\gamma=\gamma-{dg\over g}$$ since we verify that $\displaystyle A= -{\partial_x(a_0)\over a_0}-\partial_y\big({a_1\over a_0}\big)+{a_1\over a_0}\,l_1-2l_2$ and $\displaystyle B= -{\partial_y(a_0)\over a_0}-l_1$. Hence the $2$-form $k=d\gamma$ depends {\it only} on the class of the differential equation defined by $F(x,y,y')=0$. 
\medskip
\noindent
{\bf Remark 3.}\ \ The previous 1-form $\gamma$ such that $k=d\gamma$ is ``normalized'' through the case of a ${\cal L}_C(3)\subset {\bb P}^2$ associated by duality with a reduced cubic $C\subset \check{\bb P}^2$, as in Remark 2. Indeed, we get $\gamma=0$ in this case. This result can be proved by using Abel's addition theorem and is also a consequence of the relation $(*)$. Agafonov introduced an another 1-form $\gamma_{Ag}$ for a 3-web ${\cal W}(3)$ presented by $F$, such that $d\gamma_{Ag}=k$ ({\it cf.} for example [Ag-2012]). A calculus proved that the previous $\gamma$ and $\gamma_{Ag}$ are related by the following formula : $\displaystyle\gamma=\gamma_{Ag}-{1\over 2}\cdot{d\Delta\over \Delta}$ where here $\displaystyle {d\Delta\over \Delta}={du\over u}+\sum_qm_q{d\Delta_q\over \Delta_q}\cdot$ 
\medskip
\noindent
$\bullet$ {\it Dimension of the local system $\cal S$ for $d=3$}
\smallskip
For $d=3$ the system $(\cal S)$ which gives the symmetries of ${\cal W}(3)$ is the following :
$$({\cal S})\quad\left\{\matrix{
-\,\partial_x(\beta)&+&g_3\alpha+h_3\beta=0\hfill\cr
\partial_x(\alpha)-\partial_y(\beta)&+&g_2\alpha+h_2\beta=0\hfill\cr
\partial_y(\alpha)&+&g_1\alpha+h_1\beta=0.\hfill\cr
}\right.$$
Here the matrix of symmetries of the presentation $F$ is $(S_{ij})=\pmatrix{
g_3&h_3\cr
g_2&h_2\cr
g_1&h_1\cr
 }$ and there are no compatibility conditions, contrary to the case $d\geq 4$. 
\smallskip
By using methods ``à la Cartan-Spencer'' detailed in [H\'e-2004], we incarnate the local system $\cal S$ as the horizontal sections of a non necessary integrable connection $(E,\nabla)$ with rank 3, which in fact is meromorphic on the discriminant locus $|\Delta|$. This connection $(E,\nabla)$, that is equivalently a free ${\cal O}$-module of rank 3 endowed with a connection $\nabla:E\longrightarrow \Omega^1\otimes_{\cal O}E$ is called the {\it connection of symmetries} associated with $F$. It is constructed by using Cartan prolongations of the linear differential system $({\cal S})$ and the first Spencer complex on suitable jets. Here $\nabla$ is represented in an {\it adapted basis} $(e_{\ell})$ of $E$ by
$$\Gamma=\pmatrix{
-h_3dx+g_1dy&\xi_{11}dx+\xi_{12}dy&\xi_{21}dx+\xi_{22}dy\cr
-dx&g_2dx+g_1dy&h_1dy \cr
-dy&-g_3dx&-h_3dx-h_2dy\cr
}$$
with explicitly
$$\eqalign{
\xi_{11}&=(g_1+h_2)g_3-\partial_y(g_3)\cr
\xi_{12}&=g_3h_1+\partial_x(g_1)-\partial_y(g_2)\cr
\xi_{21}&=g_3h_1+\partial_x(h_2)-\partial_y(h_3)\cr
\xi_{22}&=(g_2+h_3)h_1+\partial_x(h_1)\cr
}$$
where $E:=\hbox{Ker}\, j_1=(e_1,e_2,e_3)\subset J_2({\cal O}^2)$, as jets. Without all details, we have 
$$e_1=\pmatrix{
0&-1&0&g_2-h_3&g_1&h_1 \cr
0&0&1&g_3&h_3&h_2-g_1 \cr
},$$
$$e_2=\pmatrix{
-1&g_2&g_1&\cdots \cr
0&g_3&0&\cdots \cr
}\quad \hbox{and}\quad e_3=\pmatrix{
0&0&h_1&\cdots \cr
1&h_3&h_2&\cdots \cr
}.$$ Here $j_1: J_2({\cal O}^2)\longrightarrow J_1({\cal O}^3)$ is the prolongation of the initial $j_0: J_1({\cal O}^2)\longrightarrow {\cal O}^3$ corresponding to the linear differential operator $\rho: {\cal O}^2\longrightarrow {\cal O}^3$ associated with the system $({\cal S})$, that is $j_0\pmatrix{z_1&p_1&q_1\cr
z_2&p_2&q_2\cr}
=\pmatrix{
-p_2+g_3z_1+h_3z_2\hfill\cr
p_1-q_2+g_2z_1+h_2z_2\hfill\cr
q_1+g_1z_1+h_1z_2\hfill\cr
}$ with Monge's notations.
A {\it horizontal section} $(f)=\,^t(f_1,f_2,f_3)\in \hbox{Ker}\,\nabla$ verify $df+\Gamma\,.\,f=0$ in matrix form where 
$$\left\{\matrix{
f_1=\partial_x(\alpha)+g_2\alpha=\partial_y(\beta)-h_2\beta\cr
f_2=\alpha\hfill\cr
f_3=\beta\hfill\cr
}\right.$$ corresponds by construction to $X=\alpha\partial_x+\beta\partial_y\in {\cal S}$. Morever in the basis $(e_{\ell})$, the curvature $K$ of $(E,\nabla)$ has a {\it convenient} curvature matrix
$d\Gamma+\Gamma\wedge \Gamma=\pmatrix{
k_1&k_2&k_3\cr
0&0&0 \cr
0&0&0 \cr
}dx\wedge dy$
where we have explicitly 
$$\eqalign{
k_1&=2\partial_x(g_1)-\partial_x(h_2)-\partial_y(g_2)+2\partial_y(h_3)\cr
k_2&=\partial_x(\xi_{12})-\partial_y(\xi_{11})-(g_2+h_3)\xi_{12}+g_3\xi_{22}\cr
&=\partial_x^2(g_1)-\partial_x\partial_y(g_2)+\partial_y^2(g_3)-(g_2+h_3)\big(\partial_x(g_1)-\partial_y(g_2)\big)\cr
&\quad\  -\partial_y\big((g_1+h_2)g_3\big)+\partial_x(g_3h_1)+g_3\partial_x(h_1)\cr
k_3&=\partial_x(\xi_{22})-\partial_y(\xi_{21})-(g_1+h_2)\xi_{21}+h_1\xi_{11}\cr
&=\partial_x^2(h_1)-\partial_x\partial_y(h_2)+\partial_y^2(h_3)-(g_1+h_2)\big(\partial_x(h_2)-\partial_y(h_3)\big)\cr
&\quad\ +\partial_x\big((g_2+h_3)h_1\big)-\partial_y(g_3h_1)-h_1\partial_y(g_3).\cr
}$$

By using the explicit expressions given above, we verify by a direct calculus that we have
$\displaystyle h_2=2g_1+2B+{\partial_y(\Delta)\over 2\Delta}$ and $\displaystyle h_3={g_2\over 2}+A+{\partial_x(\Delta)\over 4\Delta}\cdot$ Hence in the adapted basis $(e_{\ell})$, we obtain the {\it trace relation}
$$
\hbox{tr}(\Gamma)=-2\gamma-{d\Delta\over 2\Delta}
$$where
$\displaystyle \gamma=Adx+Bdy$. In particular, we get $\hbox{tr}(K)=d\hbox{tr}(\Gamma)=k_1dx\wedge dy=-2k$
where  $k=d\gamma$ is the curvature of the 3-web ${\cal W}(3)$ presented by $F$.
\smallskip
\noindent
{\bf Theorem 2.}\ \    {\it Let ${\cal W}(3)$ be an implicit planar $3$-web presented by $F$ with discriminant locus $|\Delta|$ and local system of symmetries $\cal S$. The following properties hold:
\smallskip
$a)$ There exists a connection $(E,\nabla)$ with rank $3$, meromorphic on $|\Delta|$, with $\hbox{Ker}\,\nabla={\cal S}$ outside $|\Delta|$ and such that its curvature $K$ vanishes if, and only if, $\dim {\cal S}=3$;
\smallskip
$b)$ There exists an explicit $3\times 3$-matrix $(k_{m\ell})$ which depends only on the class of the differential equations defined by $F(x,y,y')=0$ such that} $\dim {\cal S}=\hbox{corank}\,(k_{m\ell})$. {\it In particular $\dim {\cal S}\geq 1\Longleftrightarrow \det(k_{m\ell})=0$, {\it hence a general ${\cal W}(3)$ has no symmetry. Moreover} $\dim {\cal S}=3\Longleftrightarrow (k_{m\ell})=0\Longleftrightarrow {\cal W}(3)$ is hexagonal}.
\smallskip
\noindent
{\it Proof}. a) With the help of Cauchy-Kowalevski theorem, this part summarizes the previous results on the local system $\cal S$.

\noindent
b) From the above relation on trace, ${\cal W}(3)$ is hexagonal if $\dim {\cal S}=3$. Moreover with the expressions 
$$\eqalign{
\widetilde {\Delta_x}&=(9a_0a_3^2-7a_1a_2a_3+2a_2^3)\partial_x(a_0)+(3a_0a_2a_3-a_1a_2^2+2a_1^2a_3)\partial_x(a_1)\cr
&\quad\quad+(-3a_0a_1a_3-2a_0a_2^2+a_1^2a_2)\partial_x(a_2)+(7a_0a_1a_2-9a_0^2a_3-2a_1^3)\partial_x(a_3)\cr
\widetilde {\Delta_y}&=(9a_0a_3^2-7a_1a_2a_3+2a_2^3)\partial_y(a_0)+(3a_0a_2a_3-a_1a_2^2+2a_1^2a_3)\partial_y(a_1)\cr
&\quad\quad+(-3a_0a_1a_3-2a_0a_2^2+a_1^2a_2)\partial_y(a_2)+(7a_0a_1a_2-9a_0^2a_3-2a_1^3)\partial_y(a_3),\cr
}$$
we verify that we have
$\displaystyle k_2={1\over 2}\,\partial_x(k_1)-{\widetilde {\Delta_x}\over 2\Delta}\,k_1$ and $\displaystyle k_3={1\over 2}\,\partial_y(k_1)+{\widetilde {\Delta_y}\over 2\Delta}\,k_1$. Hence $\dim {\cal S}=3$ if, and only if, ${\cal W}(3)$ is hexagonal.  

\noindent
We also give in this part a process to describe $\dim {\cal S}$ as the corank of an effective matrix $(k_{m\ell}):{\cal O}^3\longrightarrow {\cal O}^3$ by a result proved by Olivier Ripoll in his Thesis [R-2005] ({\it cf.} also [R-2005bis]). It uses basically that ${\cal S}$ is a local system and the curvature matrix above coupled with Nakayama's lemma. The result goes as follows. For $(f)=\,^t(f_1,f_2,f_3)\in \hbox{Ker}\,\nabla$ we have $K\,.\,f=0$, that is only $k_1f_1+k_2f_2+k_3f_3=0$. Therefore, by using $\partial_x$ and $\partial_y$ with substitutions we get a matrix 
$$(k_{m\ell}):=\pmatrix{
k_1&k_2&k_3\cr 
k_{21}&k_{22}&k_{23}\cr
k_{31}&k_{32}&k_{33}\cr
}$$with two other rows given by 
$$\left\{\matrix{
k_{21}=\partial_x(k_1)+h_3k_1+k_2\hfill\cr
k_{22}=\partial_x(k_2)-\xi_{11}k_1-g_2k_2+g_3k_3\hfill\cr
k_{23}=\partial_x(k_3)-\xi_{21}k_1+h_3k_3\hfill \cr
}\right.\ \ \ \hbox{and}\ \ \ \left\{\matrix{
k_{31}=\partial_y(k_1)-g_1k_1+k_3\hfill\cr
k_{32}=\partial_y(k_2)-\xi_{12}k_1-g_1k_2\hfill\cr
k_{33}=\partial_y(k_3)-\xi_{22}k_1-h_1k_2+h_2k_3\hfill\cr
}\right.$$
such that 
$$\dim {\cal S}=\hbox{corank}\,(k_{m\ell}).$$
The presentations $F$ and $g\,.\,F$ of the ${\cal W}(3)$ at stake give the same line $(k_1,k_2,k_3)$ for $g\in {\cal O}^*$ or the same convenient curvature matrix $d\Gamma+\Gamma\wedge \Gamma$ of $(E,\nabla)$ since essentially it is true for $G$ and $H$. Hence by construction the matrix $(k_{m\ell})$ depends only on the class of the differential equations defined by $F(x,y,y')=0$. Which ends the proof of b). \cqfd 
\medskip
\noindent
{\bf Remark 4.}\  The previous result can be view as a complement of Theorem 0.  We recover in part b) a result obtained by \'Elie Cartan in the generic case in [C-1908]. The above matrix $(k_{m\ell})$ provides a new invariant associated with $F(x,y,y')=0$ and an effective method to find $\dim {\cal S}$.
\medskip
\noindent
{\bf Remark 5.}\ The symmetries $\cal S$ of a nonsingular hexagonal planar $3$-web ${\cal W}(3)$ is a local system of Lie algebras isomorphic to ${\goth b}=\{\partial_x,\partial_y,x\partial_x+y\partial_y\}$, that is of type V in the Bianchi classification from Theorem 0 and Example 1. In particular $\goth b$ is solvable with the derived series $\{0\}\subset [{\goth b},{\goth b}]=\{\partial_x,\partial_y\}\subset {\goth b}$ where $\dim [{\goth b},{\goth b}]=2$. Elements in $\goth b$ can be view in ${\goth gl}(3,{\bb C})$ as a matrix $\hbox{ad}\big(a\partial_x+b\partial_y+c(x\partial_x+y\partial_y)\big)=\pmatrix{
-c&0&a\cr
0&-c&b\cr
0&0&0\cr}$. The Lie group of $\goth b$ is $B=\{\pmatrix{
e^t&0&u\cr
0&e^t&v\cr 
0&0&1\cr
}\,;\,(t,u,v)\in {\bb C}^3\}\subset \hbox{GL}(3,{\bb C})$ and can be also view in the affine group of ${\bb C}^2$ generated by $\pmatrix{x\cr
y\cr}\longmapsto M\pmatrix{x\cr
y\cr}+\pmatrix{u\cr
v\cr}$ where $M\in \hbox{GL}(2,{\bb C})$ and $(u,v)\in {\bb C}^2$. These observations are exploited in a work in progress on monodromy properties of $(E,\nabla)$.
\medskip
\noindent
{\bf Detailed examples} (continued){\bf.} 
\smallskip
\noindent
{\bf2.} {\it The $3$-web formed by the pencils of lines through $3$ generic points}.
\smallskip
It is the basic singular hexagonal web $\displaystyle {\cal W}(x,y,{y\over x})\subset {\bb P}^2$, that is with slopes $\infty$, $0$ and $\displaystyle {y\over x}$ which is related to the relation with 3 terms given by the logarithm through the equality $\displaystyle x\cdot{1\over y}\cdot{y\over x}=1$. In the implicit setting, this web corresponds to 
 $$\displaystyle {\cal W}(y-x,y+x,{y\over x})={\cal L}_C(3)\subset {\bb P}^2$$
with slopes $\pm 1$ and $\displaystyle {y\over x}$ or equivalently with $C\subset \check {\bb P}^2$ locally defined by $P(q,p)=(p^2-1)q=0$, hence $F=(p^2-1)(y-px)$. Here we have $\Delta=4(x^2-y^2)^2$ with $\delta=x^2-y^2$ and the singular locus of this ${\cal L}_C(3)$ is globally defined in homogeneous coordinates by $X_0(X_1-X_2)(X_1+X_2)=0$. Moreover $\gamma=0$, hence $\dim {\cal S}=3$. The associated matrix of symmetries is  
$\displaystyle (S_{ij})={1\over \delta}\pmatrix{
-y&x\cr
0&0\cr
y&-x\cr
}$. 
We verify that a basis of $\cal S$ is given by 2 polynomial vector fields $x\partial_x+y\partial_y$ and $y\partial_x+x\partial_y$ completed by the multivalued  vector field  
$$\displaystyle \Big(x\log(y^2-x^2)+y\log{y+x\over y-x}\Big)\partial_x+\Big(y\log(y^2-x^2)+x\log{y+x\over y-x}\Big)\partial_y.$$

\noindent
{\bf 3.} {\it A family of examples ``à la Zariski'' with} $\dim {\cal S}=3$.
\smallskip
For $(m,n)\in {\bb N}^2$, we consider the family of $3$-webs presented by $F=p^3+x^my^n$. We have  $\Delta=-27x^{2m}y^{2n}$ with  
$\displaystyle\gamma= -{2m\over 3x}\,dx-{n\over 3y}\,dy$ and $k=d\gamma=0$. Moreover we have  $\displaystyle P_0= {n\over 3y}\,\cdot \,p^2+ {m\over 3x}\,\cdot\, p$ and its associated matrix of symmetries is
$(S_{ij})=\pmatrix{
0&0\cr
\cr
\displaystyle{m\over 3x}&\displaystyle{n\over 3y}\cr
\cr
0&0\cr
}.$ With the previous notations, we get $\Gamma=\pmatrix{
0&0&0\cr
\cr
-dx&\displaystyle{m\over 3x}dx&0\cr
\cr
-dy&0&\displaystyle-{n\over 3y}dy\cr
}$ with $\displaystyle\hbox{tr}(\Gamma)=-2\gamma -{d\Delta\over 2\Delta}$ and $K=0$. Moreover, we have the following partially multivalued basis for $n\not=3$ 
$${\cal S}=\{x^{-{m\over 3}}\partial_x,\,y^{n\over 3}\partial_y,\,(3-n)x\partial_x+(3+m)y\partial_y\}$$
and ${\cal S}=\{x^{-{m\over 3}}\partial_x,\,y\partial_y,\,\displaystyle{3x\over m+3}\,\partial_x+y\log y\,\partial_y\}$ for $n=3$.
\medskip\noindent
{\bf 4.} {\it An example which is equivalent to the normal form in case $\dim {\cal S}=1$ given by  \'Elie Cartan in} [C-1908].
\smallskip
We consider the $3$-web presented by  $F=(p^2-1)(p-u)$ with $u:=u(x)$. We have $\Delta=4(u^2-1)^2$ with $\displaystyle k=d\gamma={u''-u^2u''+2u(u')^2\over(u^2-1)^2}\not=0$ for a generic $u$. It associated matrix of symmetries is
$\displaystyle (S_{ij})={u'\over u^2-1}\pmatrix{
-1&0\cr
0&0\cr
1&0\cr
}$. We have $\det(k_{m\ell})=0$ and we verify directly that ${\cal S}=\{\partial_y\}$.
\medskip\noindent
{\bf 5.} {\it An other example with} $\dim {\cal S}=1$.
\smallskip
We consider the $3$-web presented by $F=p^3+x\,.\,p+y$. We have $\Delta=-4x^3-27y^2$ with 
$\displaystyle \gamma= -{8x^2dx+18ydy\over 4x^3+27y^2}$ and $k=d\gamma\not=0$. It associated matrix of symmetries is
$\displaystyle (S_{ij})={1\over 4x^3+27y^2}\pmatrix{
6xy&-4x^2\cr
2x^2&9y\cr
9y&-6x\cr
}.$ With the previous notations, we verify that $\det(k_{m\ell})=0$ and we get ${\cal S}=\{{\goth X}_1=2x\partial_x+3y\partial_y\}$ which can be viewed in $\hbox{Der}(\log|\Delta|)=\{{\goth X}_1,{\goth X}_2=9y\partial_x-2x^2\partial_y\}$ from Proposition 2.
\bigskip
{\bf 4. Symmetries for $d$-planar webs with $d \geq 3$ and rank problems}
\smallskip
Let ${\cal W}(d)$ be a $d$-web in $({\bb C}^2,0)$ implicitly presented by $F\in {\cal O}[p]$ with $d\geq 3$. For the $3$-upper part of  $({\cal S})$, that is the differential system at stake and from the previous paragraph, we get a connection denote again by $(E,\nabla)$. Its horizontal sections describe, up to the compatibility conditions above, the vector fields $X=\alpha\partial_x+\beta\partial_y\in {\cal S}$. In an adapted basis, this connection $(E,\nabla)$ is represented by $$\Gamma=\pmatrix{
-h_ddx+g_{d-2}dy&\xi_{11}dx+\xi_{12}dy&\xi_{21}dx+\xi_{22}dy\cr
-dx&g_{d-1}dx+g_{d-2}dy&h_{d-2}dy \cr
-dy&-g_ddx&-h_ddx-h_{d-1}dy\cr
}$$
where $\xi_{11}=(g_{d-2}+h_{d-1})g_d-\partial_y(g_d)$ and so on, with a curvature $\pmatrix{
k_1&k_2&k_3\cr
0&0&0 \cr
0&0&0 \cr
}dx\wedge dy$
where we have explicitly $k_1=2\partial_x(g_{d-2})-\partial_x(h_{d-1})-\partial_y(g_{d-1})+2\partial_y(h_d)$ and so on. Moreover $\Gamma$ gives rise to a matrix
$(k_{m\ell})=\pmatrix{
k_1&k_2&k_3\cr 
k_{21}&k_{22}&k_{23}\cr
k_{31}&k_{32}&k_{33}\cr
}$ such that $\hbox{corank}\,(k_{m\ell})$ is the ${\bb C}$-dimension of the solution of the $3$-upper part of the system $({\cal S})$.
\smallskip
 In general, the $d-3$ compatibility conditions which appear in the system $({\cal S})$ can be written 
$${g_1\over h_1}={g_2\over h_2}=\cdots={g_{d-3}\over h_{d-3}}=-{\beta\over \alpha}\cdot$$
Therefore in order to have $\dim {\cal S}\geq 1$, all the previous conditions induce specific {\it constraints} on the ordered pair $(G,H)$ associated with the implicit presentation $F$.
\medskip 
\noindent
$\bullet$ {\it Weighted Euler symmetry}
\smallskip
For a ${\cal W}(d)$ with $\dim {\cal S}\geq 1$, it may happen as the previous examples show that there exists a weighted Euler symmetry $X\in {\cal S}$, that is a vector field 
$$X=w_xx\partial_x+w_yy\partial_y\ \ \hbox{with weight} \ \ (w_x,w_y)\in {\bb C}^2-\{0\}.$$
Using the differential system $(\cal S)$, such a weighted Euler symmetry $X$ exists for a  ${\cal W}(d)$ presented by $F$ if, and only if, the system
$$({wE})\quad\left\{\matrix{
w_xxg_d+w_yyh_d&=&0\hfill\cr
w_xxg_{d-1}+w_yyh_{d-1}&=&w_y-w_x\hfill\cr
w_xxg_{d-2}+w_yyh_{d-2}&=&0\hfill\cr
&\vdots&\hfill\cr
w_xxg_1+w_yyh_1&=&0\hfill\cr
}\right.$$
has a nonzero solution $(w_x,w_y)\in {\bb C}^2$. In this case, we have necessarily 
$$\hbox{rank} \pmatrix{g_1&\ldots&g_{d-2}&g_d\cr
h_1&\ldots&h_{d-2}&h_d\cr
}\leq 1$$and for example we get a {\it radial symmetry} $X_r:=x\partial_x+y\partial_y$ if, and only if, we have $$xG+yH=0.$$ 

From previous observations, the symmetries $\cal S$ contains at most two $\bb C$-independent weighted Euler symmetries. For a sharp example suggested by Mar\'\i n, take $\displaystyle F=\prod_{i=1}^d(xp-\lambda_iy)$ with complex numbers $\lambda_i\not=\lambda_j$. Indeed, we verify here that ${\cal S}=\{x\partial_x,y\partial_y,x\log x\,\partial_x+y\log y \,\partial_y\}$ since $\displaystyle G=-{p\over x}$ and  $\displaystyle H={p\over y}\cdot$

\smallskip

We recall that a W-curve is invariant under a sub-group of $\hbox{PGL}(3,\bb C)$. It is proved in [KL-1871] that these not necessarily algebraic curves $C\subset \check{\bb P}^2$, also called anharmonic by Georges-Henri Halphen, are generated by suitable products of homogeneous equations $$X_0^{\rho_0}X_1^{\rho_1}X_2^{\rho_2}=\lambda$$ where $\lambda\in {\bb C}$ and $\rho_0+\rho_1+\rho_2=0$ with $\rho_j\in {\bb C}$. For example, let $C\subset\check{\bb P}^2$ be the algebraic W-curve of degree $d\geq 3$ with affine equation 
$$P(q,p)=p^d-\lambda q^a=0\ \ \hbox{where}\ \ 0\leq a\leq d\ \ \hbox{and}\ \ \lambda\in {\bb C}^*.$$ 
It gives rise to the $d$-web ${\cal L}_{C}(d)\subset {\bb P}^2$ endowed with a nonzero $\cal S$ which contains the weighted Euler symmetry
$$X=(d-a)\,.\,x\partial_x+d\,.\,y\partial_y.$$
Indeed, we have ${\cal S}=\{\partial_x,\,\partial_y,\,x\partial_x+y\partial_y\}$ for $a=0$ from Example 1. Otherwise, the dual curve $\check C$ is defined by a factor $\check P=\mu\lambda x^d+\nu y^{d-a}$ of the $p$-discriminant of $F(x,y,p)=P(y-px,p)$. Hence $d(d-a)\check P=X(\check P)$ for the vector field $X$ above and $X\in {\cal S}$ by construction of ${\cal L}_{C}(d)$. Moreover we may verify, for $a=d$ and at least $3\leq d\leq 6$ that ${\cal S}=\{y\partial_x,\,y\partial_y,\,y(x\partial_x+y\partial_y)\}$, and for $0<a<d$ and at least $4\leq d\leq 6$ that $\dim {\cal S}=1$. 
\medskip 
\noindent
$\bullet$ {\it Rank and remarkable planar webs}
\smallskip
Let ${\cal W}(d)={\cal W}(F_1,\ldots,F_d)$ be a nonsingular $d$-web in $({\bb C}^2,0)$ with $d\geq 3$, implicitly presented by $F$. It is defined by the 1-forms $\omega_i=dy-p_idx$ or the vector fields $X_i=\partial_x+p_i\partial_y$ where $F(x,y,p_i)=0$. A relation $\xi_1(F_1)dF_1+\cdots+\xi_d(F_d)dF_d=0$ with $\xi_i\in {\bb C}\{t\}$ is called an {\it abelian relation} of ${\cal W}(d)$. In fact these relations viewed as special $d$-uple $\big(\xi_1(F_1),\ldots,\xi_d(F_d)\big)\in {\cal O}^d$ form a local system ${\cal A}(d)$ such its {\it rank}, noted $\hbox{rank}\,{\cal W}(d)$ depends only on ${\cal W}(d)$ and is bounded by $\pi_d={1\over 2}(d-1)(d-2)$. With the same notation, the rank of a planar $d$-web ${\cal W}(d)$ is defined by using its generic nonsingular web. The previous bound is optimal. Indeed, let ${\cal L}_C(d)\subset {\bb P}^2$ be the $d$-web associated by duality with a reduced algebraic curve $C\subset \check {\bb P}^2$ with degree $d$ which is introduced in Remark 2. All the leaves of ${\cal L}_C(d)$ are straight lines and from Abel's addition theorem, we have $\hbox{rank}\,{\cal L}_C(d)=\dim_{\bb C}H^0(C,\omega_c)=\pi_d$.

For any $d$ with necessary $d\geq 5$, Mar\'\i n, Pereira and Pirio prove in [MPP-2006] that there exist {\it remarkable} webs ${\cal E}(d)$ in $({\bb C}^2,0)$. Such a $d$-web ${\cal E}(d)$ has maximum rank $\pi_d$ and is not linearisable. In particular, by a converse of Abel's addition theorem, ${\cal E}(d)$ does not come from an algebraic curve in $\check {\bb P}^2$. The classification of these ${\cal E}(d)$, even for $d=5$ is widely open.
\smallskip
If $Z=\alpha\partial_x+\beta\partial_y$ is transverse to a $d$-web ${\cal W}(d)$ in $({\bb C}^2,0)$ presented by $F$, that is  $\beta-\alpha p_i\not=0$ for $1\leq i\leq d$ where $F(x,y,p_i)=0$ , we may consider the planar $(d+1)$-web ${\cal W}(d)\sqcup Z$ presented by 
$$F_Z:=(\alpha\,.\,p-\beta)\,.\,F$$
if in addition $\alpha\not=0$. In this case ${\cal W}(d)$ appears as a sub-$d$-web of ${\cal W}(d)\sqcup Z$.
\smallskip
If a $d$-web ${\cal W}(d)$ is 
endowed with a {\it transverse symmetry} $X=\alpha\partial_x+\beta\partial_y$. Then from Theorem 1, $\displaystyle u_i(x,y):=\int_z^{(x,y)}{\omega_i\over i_X(\omega_i)}$ is generically well defined locally since $\displaystyle {\omega_i\over i_X(\omega_i)}$ is closed. Moreover $\displaystyle du_i={\omega_i\over i_X(\omega_i)}={dy-p_idx\over \beta-\alpha p_i}$ and $u_i(z)=0$ with ${\cal L}_X(u_i)=i_X(du_i)=X(u_i)=1$. Hence, if  $\displaystyle \sum_{i=1}^d\xi_i(u_i)du_i=0$ we have $\displaystyle {\cal L}_X\big(\sum_{i=1}^d\xi_i(u_i)du_i\big)=\sum_{i=1}^d\xi_i'(u_i)du_i=0$. Therefore ${\cal L}_X$ gives rise to a $\bb C$-linear map 
$${\cal L}_X: {\cal A}(d)\longrightarrow  {\cal A}(d)$$defined by $\big(\xi_i(u_i)\big)_i\longmapsto \big(\xi_i'(u_i)\big)_i$ on the $\bb C$-vector space ${\cal A}(d)$ of abelian relations of ${\cal W}(d)$. 

An efficient method introduced in [MPP-2006] to construct families of remarkable webs ${\cal E}(d)$, from algebraic W-curves $C\subset \check{\bb P}^2$, is based on the following result which used the linear map ${\cal L}_X$ through their eigenvalues and eigenspaces:
\medskip
\noindent
{\bf Theorem MPP.} \ {\it Let ${\cal W}(d)$ be a planar $d$-web with $d\geq 3$ which admits a transverse symmetry $X=\alpha\partial_x+\beta\partial_y$ with $\alpha\not=0$. Then we have}
$$\hbox{rank}\,\big({\cal W}(d)\sqcup X\big)=\hbox{rank}\,{\cal W}(d)+d-1.$$
{\it In particular ${\cal W}(d)$ is of maximum rank $\pi_d={1\over 2}(d-1)(d-2)$ if, and only if, ${\cal W}(d)\sqcup X$ is of maximum rank $\pi_{d+1}$}.
\medskip
\noindent
{\bf Examples and questions.}
\smallskip
\noindent
1. The polynomial $P(q,p)=p(p^2-1)q$ gives rise to a $\displaystyle {\cal L}_C(4)={\cal W}(y,y-x,y+x,{y\over x})\raise 2pt\hbox{,}$ that is with slopes $0$, $\pm 1$ and $\displaystyle {y\over x}\cdot$ It corresponds to the linear 4-web generated by the pencils of lines through 4 distincts points such that exactly 3 of them are aligned. Here we have $F=p(p^2-1)(y-px)$, hence $R_F=-4xy^2(x-y)^2(x+y)^2$ and $\Delta=4y^2(x-y)^2(x+y)^2$. We verify that its local system of symmetries $\cal S$ is only generated by the radial symmetry $X_r$, which moreover is {\it not} transverse to the the $\displaystyle {\cal L}_C(4)$ at stake. 
\medskip 
\noindent
2. For the W-curve $C$ defined by $P(q,p)=p^4-q$ and by using the previous methods, we verify that the local system of symmetries $\cal S$ associated with ${\cal L}_C(4)$ is only generated by $X=3x\partial_x+4y\partial_y$ which moreover is transverse. Indeed, here $\Delta := -27x^4-256y^3$ and we have  
$\displaystyle (S_{ij})={1\over \Delta}\pmatrix{
36x^2y&-27x^3\cr
-9x^3&-64y^2\cr
64y^2&-48xy\cr
48xy&-36x^2\cr
}.$  For the $3$-upper part of  $({\cal S})$, we get $k_1\not=0$ and $\det(k_{m\ell})=0$. Hence there exists a unique solution $(\alpha,\beta)$ for this differential system, up to a complex number. In fact $(3x,4y)$ verify the system $({\cal S})$. Then, we may prove directly by connection methods that ${\cal L}_C(4)\sqcup X$, presented by $F_X:=(3x\,.\,p-4y)(p^4+x\,.\,p-y)$, is remarkable. Let us note that it is not linearizable since its linearization polynomial is $\displaystyle P_0=-{36x^2\over \Delta}\cdot p^4+\cdots$ and therefore $\deg P_0>3$ ({\it cf.} for example [Hé-2014] for details and references). 
\medskip
\noindent
3. {\it Around an example given by Gilles Robert}. We start with the $4$-web presented by $F=(p^2-px+y)(p^2-px+y+1)$ which is coming from $P=(p^2+q)(p^2+q+1)$. Here $\Delta=(x^2-4y-4)(x^2-4y):=\Delta_1\,.\,\Delta_2$ and we verify that its $\cal S$ is only generated by the symmetry $X=2\partial_x+x\partial_y$ which is not weighted Euler, but with $X(\Delta_j)=0$ for $1\leq j\leq 2$. Moreover, we may check that $F_X:=(2p-x)\,.\,F$ gives a remarkable 5-web ${\cal E}(5)$ by connection methods.
\medskip
\noindent
4. Let ${\cal E}(d)$ be a remarkable planar web, then necessarily $d\geq 5$ and the following {\it dichotomy} for its Lie algebra $\cal S$ of symmetries holds: {\it $\dim{\goth g}$ is $0$ or $1$}, otherwise ${\cal E}(d)$ would be parallelizable from Proposition 1.
 
The pre-Bol has no symmetry for example, hence Bol's ${\cal B}(5)$ too. More precisely, the algebraic $4$-web ${\cal GP}(4)\subset {\bb P}^2$ formed by the pencils of lines through $4$ generic points has no symmetry. Indeed, we may choose the 4 vertices $(0,0)$, $(1,0)$, $(1,1)$ and $(0,1)$, hence with slopes 
$$p_1={y\over x}\raise 2pt\hbox{,}\ p_2={y\over x-1}\raise 2pt\hbox{,}\ p_3={y-1\over x-1}\raise 2pt\hbox{,}\ p_4={y-1\over x}\cdot$$
Then we verify by a calculus that its presentation $F$ gives $\det(k_{m\ell})\not=0$ for the $3$-upper part of its associated system $(\cal S)$. Here we recall that Bol's example ${\cal B}(5):= {\cal GP}(4)\sqcup Z$ given in [B-1936] where $Z=x(1-x)(1-2y)\partial_x+y(1-y)(1-2x)\partial_y$ is related to the five-term relation of the dilogarithm and is the first ${\cal E}(5)$ discovered.
 
\smallskip
There exist several examples of ${\cal E}(5):= {\cal L}_C(4)\sqcup Z$ endowed with a symmetry $X$. It is the case of a Terracini's example initially introduced in [Te-1937] and a Pirio's example ({\it cf.} for instance [MPP-2006]) implicitly presented by $F_Z:= (yp-x)(p^4-1)$ (resp. $(xp+y)(p^4-1)$) with the radial symmetry $X_r$ coming from the W-curve with equation $p^4-1=0$. In fact from the previous results and for a meromorphic germ $\xi=\xi(x,y)$, the planar 5-web presented by $F_\xi:=(p-\xi)(p^4-1)$, hence with $\xi^4\not=1$, admits the radial symmetry $X_r$ if, and only if, we have $x\partial_x(\xi)+y\partial_y(\xi)=0$. Hence to get all the ${\cal E}(5)$ presented by $F_\xi$, from the partial differential equation viewpoint, it is sufficient to use the explicit 6 differential conditions on $\xi$ given by connection methods in [R-2005] and the additional condition  $\partial_x(\xi)+\xi\partial_y(\xi)\not=0$ by using its linearization polynomial. For example the first condition, corresponding  to the vanishing of its generalized curvature, is $(\xi^4-1)\big(\xi^2\partial_x^2(\xi)-\partial_y^2(\xi)\big)-2\xi(\xi^4+1)\partial_x(\xi)^2+4\xi^3\partial_y(\xi)^2=0$.

\smallskip
With the terminology of [H\'e-2014]), we consider the map ${\goth u}:({\bb C}^2,0)\longrightarrow \check {\bb P}^{\,\pi_d-1}$ classically associated ``\`a la Poincar\'e-Blaschke'' with a remarkable ${\cal E}(d)$ or more generally a planar Bompiani $d$-web ${\cal W}_f(d)$ associated with a special analytic map $f:({\bb C}^2,0)\longrightarrow \check {\bb P}^{\,\pi_d-1}$. In both cases, what is the `` additional property '' for the geometry of $\goth u$ or $f$, if the web at stake is endowed with a symmetry? For $d=5$, this question appears already in [MPP-2006] for $\goth u$. A particular global interest can be noted if moreover $\goth u$ or $f$ is rational, since such examples exist.
\smallskip
In order to precise these problems at least for $d=5$ and with projective differential methods initiated by Alessandro Terracini [Te-1937] and following [H\'e-2014], let $f:({\bb C}^2,0)\longrightarrow \check{\bb P}^5$ be an analytic map with maximum $2$-osculation and associated $4\times 3$-matrix $(\alpha_{ij})$. We assume that $f$ gives rise to a planar Bompiani 5-web ${\cal W}_f(\infty,0\,;F)$ presented in normal form. This means in particular that this planar $5$-web contains the canonical ``warp and weft'' $2$-web ${\cal W}(x,y)$ and all its leaves are $1$-principal curves associated with $f$. From the projective definition of $f$, we may always suppose 
$$(\alpha_{ij})=\pmatrix{\alpha_{11}&\alpha_{12}&0\cr
\alpha_{21}&\alpha_{22}&0\cr
0&\alpha_{32}&0\cr
0&\alpha_{42}&\alpha_{43}\cr
}$$with $F:=\alpha_{42}\,.\,p^3-(3\alpha_{21}-3\alpha_{32}+\alpha_{43})\,.\,p^2-(\alpha_{11}-3\alpha_{22})\,.\,p+\alpha_{12}$ such its associated polynomial $P_0=l_1\,.\,p^2+l_2\,.\,p+l_3$ presented above verifies the following 3 explicit differential conditions of the first order on the $\alpha_{ij}$: 
$$l_1=2\alpha_{21}-\alpha_{32},\ l_2=\alpha_{11}-2\alpha_{22},\ l_3=-\alpha_{12}.$$
Moreover this ${\cal W}_f(\infty,0\,;F)$ corresponds to a ${\cal E}(5)$ if, and only if, in addition {\it only} its generalized curvature verifies $k_{{\cal W}_f(\infty,0\,;F)}=0$, that is an explicit differential condition of the second order on the $\alpha_{ij}$. Now by using a previous result, ${\cal W}_f(\infty,0\,;F)$ admits for example the radial symmetry $X_r$ if, and only if, $xG+yH=0$ or $xg_i+yh_i=0$ for $1\leq i\leq 3$ since here ${\cal W}(x,y)$ contains ${\cal W}_f(\infty,0\,;F)$. These conditions corresponds to 3 explicit differential conditions of the first order on the $\alpha_{ij}$. However, even with all these previous hypothesis, the geometric description of such $f\  modulo\  \hbox{PGL}(6,{\bb C})$ seems open.        
\bigskip
{\bf 5. Singularities for symmetries and abelian relations of a planar $3$-web}
\smallskip
Let $(E,\nabla)$ be the meromorphic connection of symmetries associated with the presentation $F$ of an implicit $3$-web ${\cal W}(3)$. It has rank $3$ and it is meromorphic on the zero locus $|\Delta|$ of the $p$-discriminant $\Delta$ of $F$. Let $(\det E, \det \nabla)$ be the {\it determinant connection} associated with $(E,\nabla)$. By definition $\det E$ is the line bundle $\wedge ^3E$ and its associated connection is presented by $\hbox{tr}(\Gamma)\in\Omega^1$ in the basis $e_1\wedge e_2\wedge e_3$ if $(e_{\ell})$ is a basis of $E$. Hence $(\det E, \det \nabla)$ is meromorphic on $|\Delta|$. 
\smallskip
The meromorphic connection of symmetries $(E,\nabla)$ is {\it regular singular} along $|\Delta|$ if for any transversal morphism $u: ({\bb C},0)\longrightarrow  ({\bb C}^2,z)$ at a smooth point $z\in |\Delta|$, the pullback
connection $(u^*E,u^*\nabla)$, which is always integrable, has a regular singularity at $0$ in the
usual sense in dimension one. This means for example that after a possible meromorphic change of basis of $u^*E$, the connection $u^*\nabla$ is presented by $\displaystyle {\bb B}(t){dt\over t}$ with an analytic $3\times 3$-matrix $\bb B$ or equivalently there exists a cyclic element for $(u^*E,u^*\nabla)$ with the
following matrix presentation:
$$\pmatrix{
0&0&-\phi_{3}  \hfill\cr
1&0&-\phi_{2}   \hfill\cr                
0&1&-\phi_1 \hfill\cr
}dt$$  
such that the $\phi_{\ell}\in{\bb C}\{t\}[1/t]$ verify {\it Fuchs'
conditions}, namely
$t^{\ell}\phi_{\ell}\in{\bb C}\{t\}$ for $1\leq \ell\leq 3$. It can be verified with the help of $\xi\longmapsto \,^t(\xi,\xi',\xi'')$ that solutions of the linear differential equation 
$$\xi'''(t)+\phi_1(t)\,.\,\xi''(t)+\phi_{2}(t)\,.\,\xi'(t)+\phi_{3}(t)\,.\,\xi(t)=0$$can be view as the horizontal sections of the dual connection
associated with $(u^*E,u^*\nabla)$. In classic study of families of projective varieties parametrized by
${\bb P}^1$ for example, the corresponding linear differential equation is the so-called {\it Picard-Fuchs equation} associated with the
{\it Gauss-Manin connection} of the family at stake. We indicate also that for an exact sequence of integrable meromorphic
connections
$$0\longrightarrow (E_1,\nabla)\longrightarrow (E_2,\nabla)\longrightarrow
(E_3,\nabla)\longrightarrow 0\,,$$
then $(E_2,\nabla)$ is regular singular if, and only if, $(E_1,\nabla)$ and
$(E_3,\nabla)$ are regular singular. Moreover the connection determinant $(\det
E, \det \nabla)$ associated with $(E,\nabla)$ is regular singular if
$(E,\nabla)$ is regular singular.
 
\medskip
Let $(E_a,\nabla_a)$ be the meromorphic connection on $|\Delta|$ associated with the abelian relations ${\cal A}(3)$ of ${\cal W}(3)$ through $F$. We have $E_a=\hbox{Ker}\,p_0$ where $p_0: J_1({\cal O})\longrightarrow {\cal O}^2$ is given by $p_0(z,p,q)=(p+Az,q+Bz)$ and $\nabla_a:E_a\longrightarrow \Omega^1\otimes_{\cal O} E_a$ is presented by $\gamma=Adx+Bdy$ in the basis $\varepsilon=(1,-A,-B)$. We recall that the horizontal sections of $(E_a,\nabla_a)$ are identified with ${\cal A}(3)$. In this approach related to Abel's addition theorem, any abelian relation of ${\cal W}(3)$ is interpreted as the vanishing trace associated with $\pi$ through the $\pi_i$ of a closed 1-form $\displaystyle a(x,y)\cdot{dy-pdx\over \partial_p(F)}=a\,.\,\nu$ on the surface $S$ defined by $F$, that is such $a$ is an analytic solution of the linear differential system 
$${\cal M}(3)\quad\quad\left\{\matrix{  
\partial_x(a)+A\,a=0 \cr 
\partial_y(a)+B\,a=0.\!\!  \cr 
}\right.$$
In fact, we have a basic isomorphism  ${\cal A}(3)\buildrel\sim\over{\longrightarrow}{\goth a}_F$ where 
${\goth a}_F=\{\omega=a\,.\,\nu\,;\,d\omega=0\}$ such the element $\big(\xi_i(F_i)\big)_i\in {\cal A}(3)$ corresponds to $a\,.\,\nu$ where $\displaystyle a=F\cdot\big(\sum_{i=1}^3{\xi_i(F_i)\partial_y(F_i)\over p-p_i}\big)\cdot$
\smallskip
In the hexagonal case and with the trace relation of the previous paragraph, we obtain from the properties above that $(E_a,\nabla_a)$ {\it is regular singular along $|\Delta|$ if it is the case for $(E,\nabla)$.}
\medskip
Let ${\cal H}(3)$ be a {\it singular hexagonal} planar $3$-web presented by $F$, that is $k=d\gamma=0$ and $\Delta(0)=0$. For such a ${\cal H}(3)$, we have $\dim {\cal S}=3$ and we may consider three symmetries $X=\alpha_1\partial_x+\beta_1\partial_y$, $Y=\alpha_2\partial_x+\beta_2\partial_y$ and $Z=\alpha_3\partial_x+\beta_3\partial_y$ in ${\cal S}$. With the identification of $\bigwedge^3{\bb C}^3$ and $\bb C$, we set
$$\theta:=\det(X,Y,Z)=\left|\matrix{
\partial_y(\beta_1)-h_2\beta_1&\partial_y(\beta_2)-h_2\beta_2&\partial_y(\beta_3)-h_2\beta_3\cr
\alpha_1&\alpha_2&\alpha_3\cr
\beta_1&\beta_2&\beta_3\cr
}\right|=\alpha\wedge \beta \wedge \partial_y(\beta)$$ for the {\it determinant} of the corresponding horizontal sections of $(E,\nabla)$. It is multivalued in general. Here we have $\alpha=(\alpha_{\ell})$ and $\beta=(\beta_{\ell})$ as column vectors, and we get the following characterization: 
\smallskip
\centerline {\it ${\cal S}=\{X,Y,Z\}$ if, and only if, $\theta$ is nonzero outside the singular locus $|\Delta|$.} 
\smallskip\noindent
In this case and with the previous notations and results, $\theta$ is a horizontal section of $(\det E, \det \nabla)$, that is
$\displaystyle 0=d\theta+\hbox{tr}(\Gamma)\,.\,\theta=d\theta-\theta(2\gamma +{d\Delta \over 2\Delta})$ and we have
$\displaystyle\gamma={d\theta \over 2\theta}-{d\Delta \over 4\Delta}\cdot$ Therefore $k=d\gamma=0$ hence there exists, locally and outside $|\Delta|$, an analytic solution $a$ of the previous linear differential system ${\cal M}(3)$. In particular, the closed $1$-form $\displaystyle a\,\cdot{dy-pdx\over \partial_p(F)}$ corresponds to an abelian relation of ${\cal H}(3)$ with $\displaystyle {da\over a}={d\Delta\over 4\Delta}-{d\theta \over 2\theta}\cdot$
\medskip
For a singular hexagonal $3$-web ${\cal H}(3)$ presented by $F$ with a reduced divisor $\displaystyle \delta=\prod _q\Delta_q$ associated with the $p$-discriminant $\Delta$, it may happen that $\gamma\in \Omega^1(\log |\Delta|)$ in the sense of Kyoji Saito ({\it cf.} [S-1980]). This means here that only $\delta\gamma\in \Omega^1$, since  $\delta d\gamma=0$ by hypothesis. In this case the integrable meromorphic connection $(E_a,\nabla_a)$ (resp. $(\det E,\det\nabla)$) is regular singular along $|\Delta|$. Moreover by using a property of the $1$-form $\gamma$ recall in Paragraph 3 and the exact sequence of K. Saito-Aleksandrov, we may suppose up to a change of presentation, that we have
$$\gamma=\sum_q\hbox{res}_q[{\cal H}(3)]\cdot{d\Delta_q\over\Delta_q}$$
with {\it residues} $\hbox{res}_q[{\cal H}(3)]:=\hbox{res}_{\Delta_q}(\gamma)$ of $\gamma$, in the sense of Poincar\'e-Leray, along the irreducible components $\Delta_q$ of $|\Delta|$. These are {\it complex numbers} which depend {\it only} on the class of $F(x,y,y')=0$ ({\it cf.} [H\'e-2006]). The previous formula, also called the {\it determinant formula}, is a way to encode informations related to the singularities of the local system ${\goth a}_F$ corresponding to the abelian relations of ${\cal H}(3)$.  
\smallskip
In regard to regular singularities for symmetries, a large class of planar hexagonal $3$-webs with {\it nonpositive rational residues } $\hbox{res}_q[{\cal H}(3)]$ is expected. For $X\in {\cal S}$ and with the notations of Proposition 2, we remark that in this case we get   
$\displaystyle \gamma(X)=\sum_q\hbox{res}_q[{\cal H}(3)]\,.\,\lambda_{X,q}$ through the invariance property.
\medskip
\noindent
{\bf Example 2} (continued){\bf.} For a ${\cal L}_C(3)$ presented by $F=(p^2-1)(y-px)$, we have $\Delta=4(x^2-y^2)^2$. With the basis of $\cal S$ given above, we get
$\theta= 2(x^2-y^2)$ and $\gamma=0$. Here $a=1$ corresponds to an abelian relation  of  ${\cal L}_C(3)$, which is a general fact from Abel's addition theorem.  
\medskip
\noindent
{\bf Example 3} (continued){\bf.} For the hexagonal web presented by $F=p^3+x^my^n$, we have $\Delta=-27x^{2m}y^{2n}$ and with $\Delta_1=x$ and $\Delta_2=y$, we get with corresponding residues
$$\gamma= -{2m\over 3}\cdot {d\Delta_1\over \Delta_1}-{n\over 3}\cdot {d\Delta_2\over \Delta_2}\cdot$$ 
Here the connection of symmetries $(E,\nabla)$ is {\it logarithmic} along $|\Delta|$ since there exists a basis of $E$ such that $\nabla$ is represented by a connection matrix $\Gamma\in \Omega^1(\log |\Delta|)$. For $n\not=3$, the basis of symmetries ${\cal S}=\{x^{-{m\over 3}}\,\partial_x,\ y^{n\over 3}\,\partial_y ,\ (3-n)x\partial_x+(3+m)y\partial_y\}$ gives $\displaystyle \theta=-{(m+3)(n-3)\over 3}\,x^{-{m\over 3}} y^{n\over 3}$ and, up to a complex number, $a=x^{2m\over 3} y^{n\over 3}$ corresponds to an abelian relation. It can be noted, contrary to the other symmetries of the basis above, that only the weighted Euler symmetry $X= (3-n)\,x\partial_x+(3+m)\,y\partial_y$ gives a complex number $\gamma(X)=-2m-n+{1\over 3}mn$. In case $n=3$, it is also true as before only for $X=y\partial_y$ with $\gamma(X)=-1$.
\medskip
\noindent
{\bf Example 6} (Mar\'\i n's family with parameters $(m,\lambda)\in {\bb N}\times {\bb C}^*)${\bf.}  It is the family of $3$-webs presented by
$$F=p(\lambda xp-y^m)\big(\lambda x(1-x)p-y^m\big)=0.$$
Here $\Delta=\lambda^2x^4y^{6m}$ and $\displaystyle \gamma=-\big({\lambda\over y^m}+{m\over y}\big)dy$. Hence we get a $(m,\lambda)$-family of hexagonal $3$-webs in ${\bb P}^2$ with a pole of order $m$ for $\gamma$ and residue $-(\lambda+1)\in {\bb C}$ along $\{y=0\}$ if $m=1$. For $m\not=1$ (resp. $m=1$) a basis of $\cal S$ is available with $\theta=-\lambda x^2y^me^{-2\lambda{y^{1-m}\over 1-m}}$ (resp. $\theta=-\lambda x^2y^{-2\lambda+1}$) which corresponds to a generator $\displaystyle a=y^me^{\lambda{y^{1-m}\over 1-m}}$ (resp. $a=y^{\lambda+1}$) for the abelian relations of the corresponding $3$-web. In particular the connection of symmetries $(E,\nabla)$ is not regular singular for $m\not=1$. However it is regular singular for  $m=1$, with a nonpositive rational residue if, and only if, $\lambda\in{\bb Q}$.   
\medskip
\noindent
{\bf Remark 6.}\ \ Suppose $\Delta$ is reduced. It is for example the case for a ${\cal L}_{C}(3)\subset {\bb P}^2$ associated with a smooth cubic $C\subset \check{\bb P}^2$. Then according to the order of the poles in the initial presentation $\Gamma=(\Gamma_{ij})$ of $(E,\nabla)$ in Paragraph 3, the meromorphic change of basis $P$ defined by $(\widetilde e_1,\widetilde e_2,\widetilde e_3)=(e_1,e_2,e_3)P$ where 
$P=\pmatrix{
\delta^{-1}&0&0\cr
0&1&0\cr
0&0&1\cr
}$ gives rise to a presentation
$$\widetilde \Gamma=P^{-1}\Gamma P+P^{-1}dP=\pmatrix{\Gamma_{11}-\displaystyle {d\delta\over \delta}&\delta\,\Gamma_{12}&\delta\,\Gamma_{13}\cr
-\displaystyle {dx\over \delta}& \Gamma_{22}& \Gamma_{23}\cr
-\displaystyle {dy\over \delta}& \Gamma_{32}& \Gamma_{33}\cr
}$$ 
which has at worst simple poles on $|\Delta|$. This proves in this case that the connection $(E,\nabla)$ and therefore also  $(E_a,\nabla_a)$ is regular singular along $|\Delta|$. 
\medskip
\noindent
$\bullet$ {\it Singular hexagonal $3$-webs related to WDVV-equations}
\smallskip

From geometry of 3-dimensional Frobenius manifolds only two non-equivalent normal forms appear as corresponding associativity equation of order three for $f\in {\cal O}$. These nonlinear partial differential WDVV-equations give rise to special planar webs called {\it booklet $3$-webs} by Agafonov in [Ag-20012]. Depending only from the following two basic equations with standard notations in the Frobenius world for partial derivatives 
$$\eqalign{ 
E_1:&=f_{yyy}+f_{xxx}\,.\,f_{xyy}-f_{xxy}^2=0\cr
E_2:&=f_{xxx}\,.\,f_{yyy}-f_{xxy}\,.\,f_{xyy}-1=0, \cr
}$$ we consider the corresponding $3$-web ${\cal W}_j(3)$ with $E_j=0$ respectively presented by
$$\eqalign{ 
W_1:&=f_{xyy}\,.\,p^3+2f_{xxy}\,.\,p^2+f_{xxx}\,.\,p-1\cr
W_2:&=f_{yyy}\,.\,p^3+f_{xyy}\,.\,p^2-f_{xxy}\,.\,p-f_{xxx}.\cr
}$$ 

With the previous notations and by calculus, the following identities hold for $W_1$: 
$$\eqalign{ 
2\Delta\,.\,A+\partial_x(\Delta)&=
-4\partial_x(E_1)\big(6f_{xxy}+f_{xxx}^2\big)\cr
2\Delta\,.\,B+\partial_y(\Delta)&=
-2\partial_x(E_1)\big(9f_{xyy}+2f_{xxx}\,.\,f_{xxy}\big)\cr
}$$ and for $W_2$:
$$\eqalign{ 
2\Delta\,.\,A+\partial_x(\Delta)&=-2\partial_x(E_2)\big(9f_{xxx}\,.\,f_{yyy}-f_{xxy}\,.\,f_{xyy}\big)-4\partial_y(E_2)\big(3f_{xxx}\,.\,f_{xyy}+f_{xxy}^2\big)\cr
2\Delta\,.\,B+\partial_y(\Delta)&=-4\partial_x(E_2)\big(3f_{xxx}\,.\,f_{xyy}+f_{xxy}^2\big)-2\partial_y(E_2)\big(9f_{xxx}\,.\,f_{yyy}-f_{xxy}\,.\,f_{xyy}\big).\cr
}$$
Hence from the above identities and the irreducible decomposition of the $p$-discriminant $\displaystyle \Delta=u\,.\,\prod_q\Delta_q^{m_q}$ of the presentation $W_1$ (resp. $W_2$), we get the following result:
\medskip
\noindent
{\bf Proposition 3.}\ {\it The $3$-webs ${\cal W}_1(3)$ and ${\cal W}_2(3)$ are hexagonal with $\displaystyle \gamma=-{1\over 2}\cdot{d\Delta\over \Delta}\cdot$ In particular, in the two cases, these hexagonal $3$-webs are regular singular from abelian relations viewpoint and the corresponding residues along $\Delta_q$ are nonpositive rational numbers equal to $\displaystyle -{m_q\over 2}\cdot$}
\medskip
\noindent
{\bf Examples WDVV-1.} In addition to $E_j=0$, solutions of WDVV-equations have quasihomogeneity constraints. By using Dubrovin's normal forms found in [D-1996] we get the following examples of ${\cal W}_1(3)$ through some explicit $f\in {\cal O}$, with rational residues from above. The vector field $E=t\partial_t+[(1-q_1)x+r_1]\partial_x+[(1-q_2)y+r_2]\partial_y$, with $(q_i,r_i)\in {\bb C}^2$ where $r_i\not=0$ only if $q_i=1$, attached to the Frobenius structure at stake provides the indicated symmetry: 
\smallskip\noindent
1. $\displaystyle f_1={x^2y^2\over 4}+{y^5\over 60}$ with $\displaystyle E={1\over 4}(4t\partial_t+3x\partial_x+2y\partial_y)$ gives rise to a ${\cal W}_1(3)$ presented by $F_{f_1}=x\,.\,p^3+2y\,.\,p^2-1$ where $\Delta=-27x^2+32y^3$ and $3x\partial_x+2y\partial_y\in {\cal S}$. A variant of this hexagonal web with a cusp as discriminant locus has been given by Alcides Lins Neto and Isao Nakai. Under the so-called regularity condition, it is the unique companion of the $3$-web ${\cal L}_C(3)\subset {\bb P}^2$ associated with $P(q,p)=p^3-q=0$ ({\it cf.} [N-2014] for details and references).   
\smallskip\noindent
2. $\displaystyle f_2={x^3y\over 6}+{x^2y^3\over 6}+{y^7\over 210}$ with $\displaystyle E={1\over 6}(6t\partial_t+4x\partial_x+2y\partial_y)$ gives the presentation $F_{f_2}=2xy\,.\,p^3+2(x+y^2)\,.\,p^2+y\,.\,p-1$ where $\Delta=4(2x+y^2)(2x-3y^2)^2$ and $2x\partial_x+y\partial_y\in {\cal S}$. 
\smallskip\noindent
3. $\displaystyle f_3={x^3y^2\over 6}+{x^2y^5\over 20}+{y^{11}\over 3960}$ with $\displaystyle E={1\over 10}(10t\partial_t+6x\partial_x+2y\partial_y)$ gives the presentation $F_{f_3}=x(x+2y^3)\,.\,p^3+y(4x+y^3)\,.\,p^2+y^2\,.\,p-1$ where $\Delta=-(27x+5y^3)(x-y^3)^3$ and $3x\partial_x+y\partial_y\in {\cal S}$.
\smallskip\noindent
4. $\displaystyle f_4=-{x^4\over 24}+xe^y$ with $\displaystyle E={1\over 2}(2t\partial_t+x\partial_x+3\partial_y)$ gives the presentation $F_{f_4}=e^y\,.\,p^3-x\,.\,p-1$ where $\Delta=e^y(4x^3-27e^y)$ and $x\partial_x+3\partial_y\in {\cal S}$. 
\bigskip\bigskip
\centerline{\bf References} 
\bigskip
\leftskip=\the\parindent \parindent=0pt
  
\textindent{[Ag-2012]}
{\smc S. I. Agafonov}, {\it Flat $3$-webs via semi-simple Frobenius $3$-manifolds}, J. Geom. Phys. {\bf 62} (2012), 361-367.

\textindent{[Ag-2015]}
{\smc S. I. Agafonov}, {\it Local classification of singular hexagonal $3$-webs with holomorphic Chern connection form and infinitesimal symmetries}, Geom. Dedicata {\bf 176} (2015), 87-115. 

\textindent{[AG-2000]}
{\smc M. A. Akivis} and {\smc V. V. Goldberg}, {\it Differential Geometry of Webs}, in Handbook of Differential Geometry, Vol. 1, Ed. F. J. E. Dillen and L. C. A. Verstraelen, Elsevier Science, Amsterdam, 2000, 1-152.

\textindent{[AS-1992]}
{\smc M. A. Akivis} and {\smc A. M. Shelekhov}, {\it Geometry and Algebra of Multidimensional Three-Webs}, Kluwer Academic Publishers, Dordrecht, 1992.

\textindent{[BB-1938]}
{\smc W. Blaschke} und {\smc G. Bol}, {\it Geometrie der Gewebe}, Springer, Berlin, 
1938.

\textindent{[B-1936]}
{\smc G. Bol}, {\it \"Uber ein bemerkenswertes F\"unfgewebe in der Ebene}, Abh. 
Hamburg {\bf 11} (1936), 387-393.

\textindent{[C-1908]}
{\smc \'E. Cartan}, {\it Les sous-groupes des groupes continus de transformations}, Ann. Sci. \'Ecole Norm. Sup. {\bf 25} (1908), 57-194.

\textindent{[Ch-1982]}
{\smc S. S. Chern}, {\it Web Geometry}, Bull. Amer. Math. Soc. {\bf 6} (1982), 
1-8.

\textindent{[D-1996]}
{\smc B. Dubrovin}, {\it Geometry of $2D$ topological field theories}, Springer Lecture Notes in Math. {\bf 1620} (1996), 120-348.  
 
\textindent{[GM-1993]}
{\smc M. Granger} and {\smc P. Maisonobe}, {\it A basic course on differential 
modules}, in $\cal D$-modules coh\'erents et holonomes, Travaux en cours {\bf 
45}, Hermann, Paris, 1993, 103-168.

\textindent{[H\'e-2000]}
{\smc A. H\'enaut}, {\it Sur la courbure de Blaschke et le rang des tissus de ${\bb C}^2$}, Natural Science Report, Ochanomizu Uni. {\bf 51} (2000), 11-25.

\textindent{[H\'e-2004]}
{\smc A. H\'enaut}, {\it On planar web geometry through abelian
relations and connections}, Ann. of Math. {\bf 159} (2004), 425-445. 

\textindent{[H\'e-2006]}
{\smc A. H\'enaut}, {\it Planar web geometry through abelian relations
and singularities}, in Inspired by S. S. Chern, Nankai Tracts in
Math. {\bf 11}, Ed. P. A. Griffiths, World Scientific, Singapore,
2006, 269-295.
 
\textindent{[H\'e-2014]}
{\smc A. H\'enaut}, {\it Planar webs of maximum rank and analytic projectives surfaces}, Math. Z. {\bf 278} (2014), 1133-1152.

\textindent{[KL-1871]}
{\smc F. Klein} und {\smc S. Lie}, {\it \"Uber diejenigen ebenen Kurven, welche durch ein geschlossen-\ es System von einfach unendlich vielen vertauschbaren
linearen Transformationen in sich \"ubergehen}, Math. Ann. {\bf 4} (1871), 50-84.

\textindent{[MPP-2006]}
{\smc D. Mar\'in}, {\smc J. V. Pereira} and  {\smc L. Pirio}, {\it On planar webs with infinitesimal automorphisms}, in Inspired by S. S. Chern, Nankai Tracts in
Math. {\bf 11}, Ed. P. A. Griffiths, World Scientific, Singapore,
2006, 351-364.

\textindent{[N-2014]}
{\smc I. Nakai}, {\it Webs and singularities}, Journal {\it of}$\,$ Singularities {\bf 9} (2014),151-167. 

\textindent{[PP-2015]}
 {\smc J. V. Pereira} and  {\smc L. Pirio}, {\it An Invitation to Web Geometry}, IMPA Monographs, Vol. 2, Springer, 2015. 

\textindent{[R-2005]}
{\smc O. Ripoll}, {\it G\'eom\'etrie des tissus du plan et \'equations diff\'erentielles}, Th\`ese de doctorat, Universit\'e Bordeaux 1, 2005.
  
\textindent{[R-2005bis]}
{\smc O. Ripoll},  {\it Détermination du rang des tissus du plan et autres invariants géo\-métriques}, C. R. Acad. Sci. Paris {\bf 341} (2005), 247-252.

\textindent{[S-1980]}
{\smc K. Saito}, {\it Theory of logarithmic differential forms and logarithmic vector fields}, J. Fac. Sci. Univ. Tokyo, {\bf 27} (1980), 265-291.

\textindent{[Te-1937]}
{\smc A. Terracini}, {\it Su una possibile particolarit\`a delle linee
principali di una
superficie, Note I e II}, Rend. della R. Acc. dei Lincei {\bf 26} (1937), 84-91 \& 153-158. 
\end